\documentclass{amsart}

\usepackage{amsmath,amsfonts,amssymb}

\usepackage[vcentermath]{youngtab}

\newtheorem{theorem}{Theorem}
\newtheorem{proposition}{Proposition}
\newtheorem{lemma}{Lemma}
\newtheorem{corollary}{Corollary}
\theoremstyle{definition}
\newtheorem{remark}{Remark}

\title[Special components of Springer fibers]
{On the singularity of some special components\\of Springer fibers}
\author{Lucas Fresse}
\address{Department of Mathematics, Weizmann Institute of Science, 76100 Rehovot, Israel, {\rm lucas.fresse@weizmann.ac.il}}
\thanks{Work supported in part by Minerva grant, No. 8596/1.}


\newcommand{\courbe}[7]{\qbezier(#1,#3)(#4,#7)(#5,#7)
\qbezier(#5,#7)(#6,#7)(#2,#3)}

\begin{document}

\begin{abstract}
Let $u\in\mathrm{End}(\mathbb{C}^n)$ be nilpotent.
The variety of $u$-stable complete flags is called the Springer fiber over $u$.
Its irreducible components are parameterized by a set of standard Young tableaux.
The Richardson (resp. Bala-Carter) components of Springer fibers
correspond to the Richardson (resp. Bala-Carter)
elements of the symmetric group, through Robinson-Schensted correspondence.
Every Richardson component is isomorphic to a product of standard flag varieties.
On the contrary, the Bala-Carter components are very susceptible to be singular.
First, we characterize the singular Bala-Carter components
in terms of two minimal forbidden configurations.
Next, we introduce two new families of components,
wider than the families of Bala-Carter components and Richardson components,
and both in duality via the tableau transposition.
The components in the first family are characterized by the fact that
they have a dense orbit of special type under the action of the stabilizer of $u$, whereas
all components in the second family are iterated
fiber bundles over projective spaces.
\end{abstract}

\maketitle

\section{Introduction}

\label{outline}

Let $V$ be an $n$-dimensional $\mathbb{C}$-vector space and let $u:V\rightarrow V$ be a nilpotent endomorphism.
We denote by ${\mathcal B}$ the set of complete flags, that is,
chains of vector subspaces $F=(V_0\subset V_1\subset \ldots\subset V_n=V)$ with $\dim V_i=i$ for all $i$.
Then, ${\mathcal B}$ is an algebraic projective variety.
We define
\[{\mathcal B}_u=\{F=(V_0,\ldots,V_n)\in{\mathcal B}:u(V_i)\subset V_i\ \mbox{for all $i$}\},\]
the subset of $u$-stable flags. Then ${\mathcal B}_u$ is a closed subvariety of ${\mathcal B}$
(in general non-irreducible).
The variety ${\mathcal B}_u$ identifies with the fiber over $u$ of the Springer resolution
(cf. \cite{Springer-1}, \cite{Springer-2}), it is called Springer fiber.
Springer fibers arise in geometric representation theory, in relation with
Springer Weyl group representations.
The study of their geometry sets quite challenging problems.
Among those, we study in this article the
question  of the singularity of their irreducible components.

The singularity of the components of ${\mathcal B}_u$ has been studied by only few authors.
First it was proved that, in some simple cases depending on the Jordan form of $u$,
every component of ${\mathcal B}_u$ is smooth.
J.A. Vargas \cite{Vargas} proved this for $u$ having only one non-trivial Jordan block
(this is the so-called hook case).
F. Fung \cite{Fung} established the property for $u$ having two blocks
(the so-called two-row case), proving in addition that in this case
every component is an iterated fiber bundle of base $(\mathbb{P}^1,\ldots,\mathbb{P}^1)$.
N. Spaltenstein \cite{Spaltenstein-1} and J.A. Vargas \cite{Vargas} provided the first example of a
singular component, in ${\mathcal B}_u$,
for $u$ having four blocks of lengths $(2,2,1,1)$.
J. Pagnon and N. Ressayre \cite{Pagnon-Ressayre} constructed a family of smooth components,
adjacent to Richardson components.

The present article comes as a continuation of 
two recent joint works with A. Melnikov.
In \cite{Fresse-Melnikov-1},
we complete the picture started by F. Fung and J.A. Vargas:
we prove that every component of ${\mathcal B}_u$ is smooth in exactly four cases
depending on the Jordan form of $u$:
1) the hook case; 2) the two-row case; 3) if $u$ has three Jordan blocks, two of arbitrary length and one which is trivial;
4) if $u$ has three blocks of length 2.

In \cite{Fresse-Melnikov-2}, we provide characterizations of the singular components
of a given Springer fiber ${\mathcal B}_u$, in the particular case $u^2=0$
(called two-column case).
For this case, we prove also that the singular components are rationally singular.
The study of the singularity of components of Springer fibers is the most fruitful
in the two-column case.
In their recent paper, 
N. Perrin and E. Smirnov \cite{Perrin-Smirnov} provide information on the type of the singularities,
they prove that the two-column type components are normal and have rational singularities.

\smallskip

In the present article, we study the singularity of the components of ${\mathcal B}_u$
for $u$ general, but while concentrating on some particular families of components.
This article contains two main results.

1) The so-called Bala-Carter components are in duality towards the Richardson components.
Whereas Richardson components are always smooth (they are parabolic orbits), 
Bala-Carter components are in 
many cases singular.
As a first result, we give a necessary and sufficient condition for a Bala-Carter component to be singular.
It will follow from the criterion that, among components of Springer fibers, Bala-Carter components
are the most susceptible to be singular,
in the sense that whenever ${\mathcal B}_u$ has a singular component,
it admits one of Bala-Carter type.

2) We introduce a new family of components generalizing the Bala-Carter components, which is in duality
towards a family of components generalizing the Richardson components.
The generalized Bala-Carter components are those containing a dense orbit of a special type
under the action of the stabilizer of $u$.
Then 
we prove that every generalized Richardson component is an iterated fiber bundle over projective spaces.
In the two-column case, every component is generalized Bala-Carter,
whereas in the two-row case, every component is generalized Richardson.
Then we retrieve in particular the result due to F. Fung.

\section{Background and statement of main results}

\label{section-background-statement}

Before stating our results, 
which will be done in the subsections \ref{first-result} and \ref{second-result},
we need to set up the basic background.

\subsection{Components of Springer fibers}
The variety ${\mathcal B}_u$ is an algebraic projective variety, which is connected
but, in general, reducible.
Following R. Steinberg \cite{Steinberg} and N. Spaltenstein \cite{Spaltenstein},
${\mathcal B}_u$ is equidimensional and its
irreducible components are parameterized by a set of standard tableaux.
In this subsection, first, we present Spaltenstein's construction.
Then, we recall the definition of two special families of components of ${\mathcal B}_u$: the
Richardson and Bala-Carter components.

\subsubsection{The Jordan form $\lambda(u)$ and the Young diagram $Y(u)$}

\label{section-Yu}

Let $\lambda(u)=(\lambda_1\geq\lambda_2\geq \ldots\geq\lambda_r)$ be the sizes of the Jordan blocks of $u$.
We have $\lambda_1+\lambda_2+\ldots+\lambda_r=n$, that is, the sequence $\lambda(u)$ is a partition of $n$.
Let $Y(u)$ be the Young diagram of rows of lengths $(\lambda_1, \ldots, \lambda_r)$,
that is, $Y(u)$ is an array of $r$ left-adjusted rows with the $i$-th row containing $\lambda_i$ empty boxes.
\[\mbox{Example:}\quad
\lambda(u)=(3,2,2,1)\ \Rightarrow\ Y(u)=\yng(3,2,2,1).
\]
The diagram $Y(u)$ is a datum equivalent to the sequence $\lambda(u)$, and to the Jordan form of $u$.
Also, the particular Jordan forms invoked in section \ref{outline}
may be interpreted in terms of the diagram $Y(u)$.
In a transparent way, in the hook case (i.e., $u$ has only one non-trivial Jordan block)
the diagram has only one row of length $\geq 2$,
in the two-row case (i.e., $u$ has two blocks) the diagram has two rows,
and respectively, in the two-column case (i.e., $u^2=0$) the diagram has two columns.

Let $\lambda^*(u)=(\lambda^*_1\geq\ldots\geq\lambda^*_s)$ be the partition of $n$ conjugate of $\lambda(u)$,
that is, $\lambda^*_1,\ldots,\lambda^*_s$ are the sizes of the columns of the diagram $Y(u)$.
The dimension of ${\mathcal B}_u$ has the following expression in terms of $\lambda^*(u)$
(see \cite[\S II.5]{Spaltenstein}):
\begin{equation}
\label{relation-dimension-Bu}
\dim{\mathcal B}_u=\sum_{j=1}^{s}\frac{\lambda^*_j(\lambda^*_j-1)}{2}.
\end{equation}

\subsubsection{The irreducible components ${\mathcal K}^T\subset{\mathcal B}_u$}

\label{Spaltenstein-construction}

Recall that a standard Young tableau (in short, standard tableau)
is a numbering of $Y(u)$ with the entries
$1,\ldots,n$, such that the entries increase from left to right along the rows
and respectively
from top to bottom along the columns. For instance, here is a standard tableau
of shape $Y(u)$, where $Y(u)$ comes from the previous example:
\[T=\young(138,25,46,7).\]
A standard tableau $T$ is in fact equivalent to the datum of the maximal
chain of subdiagrams $Y_1^T\subset Y_2^T\subset\ldots\subset Y_n^T=Y(u)$,
where $Y_i^T$ is the shape of the subtableau of $T$ formed by the entries $1,\ldots,i$.

Spaltenstein's construction of the components of ${\mathcal B}_u$ relies on a partition
of ${\mathcal B}_u$ into subsets ${\mathcal B}_u^T$ associated to the $T$'s standard.
For $F=(V_0,V_1,\ldots,V_n)\in{\mathcal B}_u$,
each $V_i$ is $u$-stable and the restriction
$u_{|V_i}\in\mathrm{End}(V_i)$ is nilpotent. Its Jordan form is represented by a Young subdiagram
$Y(u_{|V_i})\subset Y(u)$. This defines a maximal chain of subdiagrams
$Y(u_{|V_1})\subset Y(u_{|V_2})\subset\ldots\subset Y(u_{|V_n})=Y(u)$.
Then we set
\[{\mathcal B}_u^T=\{F\in {\mathcal B}_u:Y(u_{|V_i})=Y_i^T\ \forall i=1,\ldots,n\}.\]
The variety ${\mathcal B}_u$ is indeed the disjoint union of the subsets ${\mathcal B}_u^T$.
Due to \cite[\S II.5]{Spaltenstein}, each subset ${\mathcal B}_u^T$ is locally closed in ${\mathcal B}_u$,
irreducible and smooth, and $\dim {\mathcal B}_u^T=\dim {\mathcal B}_u$. Therefore,
the irreducible components of ${\mathcal B}_u$ are the closures ${\mathcal K}_u^T:=\overline{{\mathcal B}_u^T}$
in the Zariski topology,
they are parameterized by the standard tableaux of shape $Y(u)$,
and we have $\dim {\mathcal K}_u^T=\dim {\mathcal B}_u$ for every $T$.

Note that, up to isomorphism,
${\mathcal K}^T_u$ depends only on the tableau $T$. That is why we will drop 
the index from the notation and
write ${\mathcal K}^T={\mathcal K}^T_u$.

\subsubsection{Bala-Carter and Richardson components}

\label{def-BC-R}

A Bala-Carter component of ${\mathcal B}_u$ 
is associated to a 
permutation of $\lambda(u)$,
that is, a sequence $\pi=(\pi_1,\ldots,\pi_r)$ which coincides with $\lambda(u)=(\lambda_1,\ldots,\lambda_r)$
up to ordering. A Richardson component 
is associated to a permutation
of the conjugate partition $\lambda^*(u)$.
We denote by $\Lambda_u$ the set of permutations of $\lambda(u)$ and
by $\Lambda^*_u$ the set of permutations of $\lambda^*(u)$.

Observe preliminarily that a flag $F=(V_0,\ldots,V_n)\in{\mathcal B}_u$
induces a nilpotent endomorphism $u_{|V_j/V_i}\in\mathrm{End}(V_j/V_i)$
for any $0\leq i<j\leq n$. Moreover, the map
$F\mapsto \mathrm{rank}\,u_{|V_j/V_i}$
is lower semi-continuous (see \cite[Lemma 2.2]{Fresse}).

Let $\pi=(\pi_1,\ldots,\pi_r)\in \Lambda_u$.
For $j\in\{0,\ldots,r\}$, set $i_j=\pi_1+\ldots+\pi_j$.
We say that $F\in{\mathcal B}_u$ is $\pi$-regular if
$u_{|V_{i_j}/V_{i_{j-1}}}$ is regular for all $j=1,\ldots,r$.
The set
\[{\mathcal U}^\mathrm{BC}_{\pi}=\{F\in{\mathcal B}_u:\mbox{$F$ is $\pi$-regular}\}\]
is then an open subset of ${\mathcal B}_u$.
Let $Z_u=\{g\in GL(V):gug^{-1}=u\}$ be the stabilizer of $u$.
In fact, the set ${\mathcal U}^\mathrm{BC}_{\pi}$ is irreducible, since it is a $Z_u$-orbit of ${\mathcal B}_u$
(see also section \ref{Jordan-orbits}).
Hence, its closure in the Zariski topology,
denoted by ${\mathcal K}^\mathrm{BC}_\pi$, is an irreducible component of ${\mathcal B}_u$,
called a Bala-Carter component (see also \cite[\S 5.10]{Carter}).

Let $\pi=(\pi_1,\ldots,\pi_s)\in \Lambda^*_u$ and,
for $j\in\{0,\ldots,s\}$, set $i_j=\pi_1+\ldots+\pi_j$.
We say that $F\in{\mathcal B}_u$ is $\pi$-trivial if
$u_{|V_{i_j}/V_{i_{j-1}}}=0$ for all $j=1,\ldots,s$.
The set
\[{\mathcal K}^\mathrm{R}_{\pi}=\{F\in{\mathcal B}_u:\mbox{$F$ is $\pi$-trivial}\}\]
is then a closed subset of ${\mathcal B}_u$. In fact, there is a unique
partial flag $W_1\subset W_2\subset\ldots\subset W_s=V$ with $\dim W_j=i_j$ and $u(W_j)\subset W_{j-1}$
for all $j$, hence, $F\in{\mathcal B}_u$ is $\pi$-trivial
if and only if $V_{i_j}=W_j$ for all $j$.
Denoting by ${\mathcal B}(W)$ the variety of complete flags of a space $W$
and letting ${\mathcal B}^{(m)}={\mathcal B}(\mathbb{C}^m)$, we obtain
a natural isomorphism
\begin{equation}
\label{iso-Richardson}
{\mathcal K}^\mathrm{R}_{\pi}\cong {\mathcal B}(W_1)\times {\mathcal B}(W_2/W_1)\times\cdots\times{\mathcal B}(W_s/W_{s-1})
\cong {\mathcal B}^{(\lambda^*_1)}\times{\mathcal B}^{(\lambda^*_2)}\times\cdots\times{\mathcal B}^{(\lambda^*_s)}
\end{equation}
(see \cite[\S 7]{Pagnon-Ressayre}).
Therefore, ${\mathcal K}^\mathrm{R}_{\pi}$ is irreducible and of same dimension as ${\mathcal B}_u$.
Hence, it is an irreducible component of ${\mathcal B}_u$, called a Richardson component.

\subsubsection{Duality between Bala-Carter and Richardson components}

\label{duality}

Let us describe the standard tableaux corresponding to Richardson and Bala-Carter components.
Let $\pi=(\pi_1,\ldots,\pi_r)\in\Lambda_u$, and let $T_\pi$ be the tableau obtained as follows:
draw an array of $r$ left-adjusted rows, with $\pi_i$ boxes in the $i$-th row
filled in with the numbers $\pi_1+\ldots+\pi_{i-1}+1,\ldots,\pi_1+\ldots+\pi_{i}$,
and then push to the up the supernumerary boxes in each row, so to obtain a standard tableau of shape $Y(u)$.
For example:
\[\pi=(2,3,1,2)\ \ \rightarrow\ \ \young(12,345,6,78)\ \ \rightarrow\ \ \young(125,34,68,7)=T_\pi.\]
Then, we can see that ${\mathcal U}_\pi^\mathrm{BC}\subset{\mathcal B}_u^{T_\pi}$ (see also section \ref{Jordan-orbits}),
so that ${\mathcal K}_\pi^\mathrm{BC}={\mathcal K}^{T_\pi}$.

Similarly, for $\pi=(\pi_1,\ldots,\pi_s)\in\Lambda^*_u$, let $T^*_\pi$ be the tableau obtained as follows:
draw an array of $s$ top-adjusted columns, with $\pi_i$ boxes in the $i$-th column
filled in with the numbers $\pi_1+\ldots+\pi_{i-1}+1,\ldots,\pi_1+\ldots+\pi_{i}$,
and push to the left the supernumerary boxes in each column, to obtain a standard tableau of shape $Y(u)$.
Then, we have the inclusion ${\mathcal B}_u^{T^*_\pi}\subset{\mathcal K}_\pi^\mathrm{R}$, and therefore,
${\mathcal K}_\pi^\mathrm{R}={\mathcal K}^{T^*_\pi}$ (see \cite[\S 7]{Pagnon-Ressayre}).

In fact, the tableau $T^*_\pi$ is obtained through the Robinson insertion algorithm
from a Richardson element $w\in\mathbf{S}_n$, whereas
$T_\pi$ is obtained from a Bala-Carter element, i.e., of the form $w_0ww_0^{-1}$, with
$w\in\mathbf{S}_n$ Richardson and $w_0:i\mapsto n-i+1$
(cf., \cite{Fulton}).

To $u$, we may associate a nilpotent $u^*\in\mathrm{End}(V)$ such that $u$, $u^*$ have conjugate
Jordan forms: $\lambda(u^*)=\lambda^*(u)$. We consider the corresponding Springer fiber
${\mathcal B}_{u^*}$.
For $T$ a standard tableau of shape $Y(u)$, denote by $T^*$ its transpose of shape $Y(u^*)$, that is,
the $i$-th column of $T^*$ coincides with the $i$-th row of $T$.
Thus, to each component ${\mathcal K}^T={\mathcal K}_u^T\subset{\mathcal B}_u$,
it corresponds the component ${\mathcal K}^{T^*}={\mathcal K}_{u^*}^{T^*}\subset{\mathcal B}_{u^*}$.
From the description above, we infer the duality property:
\[\mbox{${\mathcal K}^T$ is Bala-Carter\ \ $\Leftrightarrow$\ \ ${\mathcal K}^{T^*}$ is Richardson.}\]

\subsection{Singularity criterion for Bala-Carter components}

\label{first-result}

It follows from formula (\ref{iso-Richardson}) that the Richardson components
of ${\mathcal B}_u$ are smooth and pairwise isomorphic.
On the contrary, the Bala-Carter components of ${\mathcal B}_u$ are not pairwise isomorphic in general,
in addition there can be smooth and singular ones.
Our first result is a characterization of the singular Bala-Carter components.

We recall two examples of singular components of Springer fibers.
The first example is ${\mathcal K}^T\subset{\mathcal B}_u$ for $u$ having the Jordan form $\lambda(u)=(2,2,1,1)$,
the second example is ${\mathcal K}^S\subset{\mathcal B}_u$ for $u$ of Jordan type $\lambda(u)=(3,2,2)$.
These two singular components are associated to the tableaux
\[T=\young(13,25,4,6)\quad\mbox{and}\quad S=\young(125,34,67)\]
(see \cite{Fresse-Melnikov-1} and \cite{Vargas}).
Notice that both components are Bala-Carter, associated to the
following permutations of the Jordan block size sequences:
$(1,2,2,1)$ and $(2,3,2)$.

If $\pi=(\pi_1,\ldots,\pi_r)$ and $\rho=(\rho_1,\ldots,\rho_k)$ are two sequences of integers,
say $\pi\geq\rho$ if there are $1\leq i_1<\ldots<i_k\leq r$
such that $\pi_{i_l}\geq\rho_l$ for every $l\in\{1,\ldots,k\}$
(in particular $r\geq k$).
Then, our result can be stated:

\begin{theorem}
\label{theorem-1}
Let $\pi\in\Lambda_u$ be a permutation of the Jordan block sizes sequence $\lambda(u)$, and
let ${\mathcal K}_\pi^{\mathrm{BC}}\subset{\mathcal B}_u$ be the corresponding Bala-Carter component.
Then, ${\mathcal K}_\pi^{\mathrm{BC}}$ is singular if and only if
$\pi\geq (1,2,2,1)$ or $\pi\geq(2,3,2)$.
\end{theorem}

We derive two corollaries from this characterization.
First, as we know from \cite{Fresse-Melnikov-1} that ${\mathcal B}_u$ admits a singular component
provided that $\lambda(u)\geq (2,2,1,1)$ or $\lambda(u)\geq (3,2,2)$,
we see that in this case, we can always find $\pi\in\Lambda_u$ satisfying $\pi\geq (1,2,2,1)$ or $\pi\geq (2,3,2)$.
That is:

\begin{corollary}
Whenever the Springer fiber ${\mathcal B}_u$ admits a singular component, it admits a singular component of Bala-Carter type.
\end{corollary}

Second, we easily infer from the theorem that in many cases, all the Bala-Carter components of ${\mathcal B}_u$ are
singular:

\begin{corollary}
If $\lambda(u)\geq(2,2,2,2)$ or $\lambda(u)\geq(3,3,3)$, then every Bala-Carter component of ${\mathcal B}_u$ is singular.
\end{corollary}

It is convenient to represent a sequence $\pi=(\pi_1,\ldots,\pi_r)\in\Lambda_u$ by a graph:
each number $\pi_i$ is represented by a chain of $\pi_i-1$ arcs connecting two by two
$\pi_i$ vertices arranged along a horizontal line, then the graph
corresponding to $\pi$ is obtained by juxtaposing the chains for $\pi_1,\ldots,\pi_r$.
For instance, $(1,2,2,1)$ and $(2,3,2)$ are represented by the graphs
\begin{equation}\label{two-singular-graphs}\begin{picture}(80,10)(-16,0)
\courbe{0}{16}{0}{7}{8}{9}{8}
\courbe{32}{48}{0}{39}{40}{41}{8}
\put(-18,-3){$\bullet$}
\put(-2,-3){$\bullet$}
\put(14,-3){$\bullet$}
\put(30,-3){$\bullet$}
\put(46,-3){$\bullet$}
\put(62,-3){$\bullet$}
\end{picture}
\qquad\mbox{and}\qquad
\begin{picture}(96,10)(-16,0)
\courbe{-16}{0}{0}{-9}{-8}{-7}{8}
\courbe{16}{32}{0}{23}{24}{25}{8}
\courbe{32}{48}{0}{39}{40}{41}{8}
\courbe{64}{80}{0}{71}{72}{73}{8}
\put(-18,-3){$\bullet$}
\put(-2,-3){$\bullet$}
\put(14,-3){$\bullet$}
\put(30,-3){$\bullet$}
\put(46,-3){$\bullet$}
\put(62,-3){$\bullet$}
\put(78,-3){$\bullet$}
\end{picture}\ \ .\end{equation}
The relation $\pi\geq\rho$ is then translated in terms of an inclusion of graphs: the graph of
$\rho$ can be obtained from the one of $\pi$ by
repeating the procedure of deleting
either a whole connected component or the extremal vertex of a connected component (together with the corresponding arc).
Then, the Bala-Carter component ${\mathcal K}_\pi^{\mathrm{BC}}\subset{\mathcal B}_u$ is singular
if and only if the graph representing $\pi$ contains one of both graphs
of formula (\ref{two-singular-graphs}).

\subsection{Generalization of Bala-Carter and Richardson components}

\label{second-result}

Let $Z_u=\{g\in GL(V):gug^{-1}=u\}$ be the stabilizer of $u$.
This is a connected closed subgroup of $GL(V)$, and its natural action on flags
leaves ${\mathcal B}_u$ and every component of ${\mathcal B}_u$ invariant.
In this subsection, we define a notion of particular orbit of ${\mathcal B}_u$
under the action of $Z_u$, which we will call Jordan orbit.
We point out that
a component is Bala-Carter if and only if it contains a dense Jordan orbit of a particular
type, called standard.
Then we consider the components which contain a (not necessarily standard) Jordan orbit,
they provide a generalization of Bala-Carter components.
We state two results: the first one characterizes the generalized Bala-Carter components,
the second one says that, if ${\mathcal K}^T$ is generalized Bala-Carter,
then the component ${\mathcal K}^{T^*}$ is an iterated fiber bundle of base a sequence of projective spaces.

\subsubsection{Jordan orbits}

\label{Jordan-orbits}

Recall that $\lambda(u)=(\lambda_1,\ldots,\lambda_r)$ are the sizes of the Jordan blocks of $u$.
We denote by $\Pi_u$ the set of partitions
of $\{1,\ldots,n\}$ into $r$ subsets of cardinal $\lambda_1,\ldots,\lambda_r$.
Hence, an element $\pi\in\Pi_u$ can be written $\pi=(I_1,\ldots,I_r)$ such that
\[I_1\sqcup\ldots\sqcup I_r=\{1,\ldots,n\}\quad\mbox{and}\quad |I_k|=\lambda_k\mbox{\ \ for all $k$}.\]
We will write $I_1,\ldots,I_r\in\pi$.
A partition $\pi\in\Pi_u$ can be (and will be) identified to a
map $\pi:\{1,\ldots,n\}\rightarrow\{1,\ldots,n\}\cup\{\emptyset\}$
such that
\[
\mbox{$\pi(i)\in\{\emptyset,1,\ldots,i-1\}$,}
\quad\mbox{and}\quad
\pi(i)=\pi(j)\not=\emptyset\ \Rightarrow\ i=j.
\]
Indeed, assuming $i\in I_k$, let $\pi(i)=\emptyset$ if $i=\mathrm{min}(I_k)$, and let $\pi(i)$
be the predecessor of $i$ in $I_k$ otherwise.
Alternatively, we represent $\pi$ by a pattern (also denoted by $\pi$)
consisting of $n$ points along a horizontal line which we label
from left to right by $1,\ldots,n$, with arcs $(j,i)$ for $\pi(i)=j$.
For instance, 
$\pi=(\{1,2,5\},\{3,8\},\{6,7\},\{4\})$ is represented by
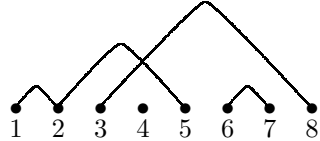
\[
\begin{picture}(112,60)(0,-12)
\courbe{0}{16}{0}{7}{8}{9}{8}
\courbe{16}{64}{0}{37}{40}{43}{24}
\courbe{32}{112}{0}{69}{72}{75}{40}
\courbe{80}{96}{0}{87}{88}{89}{8}
\put(-2,-3){$\bullet$}
\put(14,-3){$\bullet$}
\put(30,-3){$\bullet$}
\put(46,-3){$\bullet$}
\put(62,-3){$\bullet$}
\put(78,-3){$\bullet$}
\put(94,-3){$\bullet$}
\put(110,-3){$\bullet$}
\put(-2,-11){1}
\put(14,-11){2}
\put(30,-11){3}
\put(46,-11){4}
\put(62,-11){5}
\put(78,-11){6}
\put(94,-11){7}
\put(110,-11){8}
\end{picture}
\]

For $\pi\in\Pi_u$, a basis $(e_1,\ldots,e_n)$ of $V$ is said to be a $\pi$-basis
if it satisfies
\[
\left\{\begin{array}{ll} u(e_i)=e_{\pi(i)} & \mbox{if $\pi(i)\not=\emptyset$}, \\
u(e_i)=0 & \mbox{otherwise}
\end{array}\right.
\]
(this is simply a Jordan basis, numbered according to $\pi$). We denote by ${\mathcal Z}_\pi$ the set of flags which can be written
$F=(\langle e_1,\ldots,e_i\rangle)_{i=0,\ldots,n}$ for some $\pi$-basis
$(e_1,\ldots,e_n)$.
Clearly, ${\mathcal Z}_\pi\subset{\mathcal B}_u$.
Notice that the set of $\pi$-bases is an orbit of $V^n$ under the action of $Z_u$.
Therefore, ${\mathcal Z}_\pi$ is a $Z_u$-orbit of ${\mathcal B}_u$.
We call it a Jordan orbit.

To a partition $\pi\in\Pi_u$, we associate a standard tableau $T_\pi$ in the following manner.
For $i\in\{1,\ldots,n\}$, let $c_\pi(i)$ be minimal such that $\pi^{c_\pi(i)}(i)=\emptyset$.
Then, let $T_\pi$ be the unique standard tableau of shape $Y(u)$ which contains $i$ in its $c_\pi(i)$-th
column for every $i$. For instance, for $\pi$ like in the previous example, we have
\[T_\pi=\young(125,37,48,6).\]
For $\pi\in\Pi_u$,
it is easy to see that the orbit ${\mathcal Z}_\pi$ lies in the set
${\mathcal B}_u^{T_\pi}$ (cf. section \ref{Spaltenstein-construction}) and therefore, in the component ${\mathcal K}^{T_\pi}$.

\begin{remark}
In general, any $Z_u$-orbit of ${\mathcal B}_u$ is not a Jordan orbit.
Suppose for example
$\lambda(u)=(3,1)$. Let $(e,e',e'',f)$ with $u(e'')=e'$, $u(e')=e$, $u(e)=u(f)=0$ be a
Jordan basis. Then the flag $F=(0\subset\langle
e\rangle\subset\langle e,e'+f\rangle\subset\langle
e,e',f\rangle\subset V)$ does not belong to any ${\mathcal
Z}_\pi$.
However, the Jordan orbits are all the $Z_u$-orbits of ${\mathcal B}_u$
in the two-column case (see \cite[\S 2]{Fresse-Melnikov-2}).
\end{remark}

Now, let us characterize a Bala-Carter component in terms of Jordan orbits.

A partition $\pi\in\Pi_u$ is said to be standard if $\pi(i)\in\{\emptyset,i-1\}$
for all $i$
(equivalently, $\pi=(I_1,\ldots,I_r)$ where the $I_j$'s are integer intervals).
We denote by $\Pi_u^0\subset\Pi_u$ the subset of standard partitions.
The orbit ${\mathcal Z}_\pi$ associated to a standard partition $\pi$
is called a standard Jordan orbit.
Notice that there is a one-to-one correspondence between
the set $\Lambda_u$ of permutations
of the sequence $\lambda(u)$ and the set $\Pi_u^0$ of standard partitions:
an element $\pi=(\pi_1,\ldots,\pi_r)\in\Lambda_u$ can be seen as a standard partition
$(I_1,\ldots,I_r)\in\Pi_u^0$ where
for $j=1,\ldots,r$, we set $I_j=\{\pi_1+\ldots+\pi_{j-1}+1,\ldots,\pi_1+\ldots+\pi_j\}$.
The definition of the tableau $T_\pi$ for $\pi\in\Lambda_u$ in the subsection \ref{duality}
is compatible with the previous one given for $\pi\in\Pi_u$ in this subsection.

Let $\pi=(\pi_1,\ldots,\pi_r)\in\Lambda_u$ and
let us come back to the definition of the Bala-Carter component
${\mathcal K}_\pi^{\mathrm{BC}}\subset {\mathcal B}_u$.
For $j=0,\ldots,r$, write $i_j=\pi_1+\ldots+\pi_j$.
Then ${\mathcal K}_\pi^{\mathrm{BC}}$ contains as a dense subset the set
${\mathcal U}_\pi^{\mathrm{BC}}$ formed by the flags $F=(V_0,\ldots,V_n)\in{\mathcal B}_u$
such that $u_{|V_{i_j}/V_{i_{j-1}}}$ is regular for all $j$.
For such a flag, we find
$x_{j}\in V_{i_j}$ such that
$$x_j, u(x_j),\ldots,u^{\pi_j-1}(x_j)\notin V_{i_{j-1}}.$$
Due to the Jordan form of $u$, we have $V_{i_j}\subset V_{i_{j-1}}+\ker u^{\pi_j}$ for all $j$, hence we may
choose $x_j\in\ker u^{\pi_j}$.
For $i_{j-1}<i\leq i_j$, put $e_i=u^{i_j-i}(x_j)$. We necessarily have
$$V_i=V_{i_{j-1}}+\langle e_l:i_{j-1}<l\leq i\, \rangle.$$
Altogether, $e_1,\ldots,e_n$ form a basis of $V$ which satisfies
\[V_i=\langle e_1,\ldots,e_i\rangle\quad \forall i\in\{0,\ldots,n\},\]
\[
\mbox{and}\quad
\left\{\begin{array}{ll}
u(e_i)=0 & \mbox{if $i=i_j+1$ for some $j\in\{0,\ldots,r-1\}$}, \\
u(e_i)=e_{i-1} & \mbox{otherwise.}
\end{array}
\right.
\]
Thus, this is a $\pi$-basis.
Conversely, if $(e_1,\ldots,e_n)$ is a $\pi$-basis,
then the flag $(\langle e_1,\ldots,e_i\rangle)_{i=0,\ldots,n}$ lies in ${\mathcal U}_\pi^{\mathrm{BC}}$.
It results ${\mathcal U}_\pi^{\mathrm{BC}}={\mathcal Z}_\pi$,
hence ${\mathcal K}_\pi^{\mathrm{BC}}$ contains a dense standard Jordan orbit.
Using the correspondence $\Pi_u^0\cong\Lambda_u$, we therefore obtain:

\begin{proposition}
A component ${\mathcal K}^T\subset{\mathcal B}_u$ is Bala-Carter
if and only if it contains a dense standard Jordan orbit.
\end{proposition}

\subsubsection{Generalized Bala-Carter components}

\label{generalized-bc-components}

We consider the components which have the property to contain a dense Jordan orbit ${\mathcal Z}_\pi$
(not necessarily standard).
In particular, due to the above, Bala-Carter components hold this property.
Our first purpose is to characterize these components.

Let $\pi=(I_1,\ldots,I_r)\in\Pi_u$ be a partition of $\{1,\ldots,n\}$.
We say that $\pi$ has a crossing if there exist $i,j\in\{1,\ldots,n\}$ with $\emptyset<\pi(j)<\pi(i)<j<i$
(that is, there is a crossing of two arcs in the graph representing $\pi$).
We say that $I_j<I_k$ if $I_j\subset\,]\mathrm{min}(I_k),\mathrm{max}(I_k)[\,$.
Then, we define $\Pi_u^1$ as the subset of elements $\pi\in\Pi_u$ with
no crossings, and satisfying
\[I_j<I_k\ \Rightarrow\ |I_j|\geq|I_k|.\]
Note that in particular $\Pi_u^0\subset \Pi_u^1$.
The elements of $\Pi_u^1$ are characterized by the fact that their graph representation
has no crossing, and if a chain of arcs lies under another one, then the first chain is the longest.
For example, if the graph is
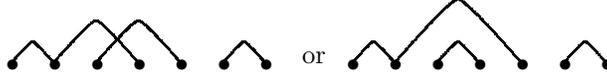
\[
\begin{picture}(100,30)(0,0)
\courbe{0}{16}{0}{7}{8}{9}{8}
\courbe{16}{48}{0}{30}{32}{34}{16}
\courbe{32}{64}{0}{46}{48}{50}{16}
\courbe{80}{96}{0}{87}{88}{89}{8}
\put(-2,-3){$\bullet$}
\put(14,-3){$\bullet$}
\put(30,-3){$\bullet$}
\put(46,-3){$\bullet$}
\put(62,-3){$\bullet$}
\put(78,-3){$\bullet$}
\put(94,-3){$\bullet$}
\end{picture}
\quad\mbox{or}\quad
\begin{picture}(100,30)(0,0)
\courbe{0}{16}{0}{7}{8}{9}{8}
\courbe{16}{64}{0}{37}{40}{43}{24}
\courbe{32}{48}{0}{39}{40}{41}{8}
\courbe{80}{96}{0}{87}{88}{89}{8}
\put(-2,-3){$\bullet$}
\put(14,-3){$\bullet$}
\put(30,-3){$\bullet$}
\put(46,-3){$\bullet$}
\put(62,-3){$\bullet$}
\put(78,-3){$\bullet$}
\put(94,-3){$\bullet$}
\end{picture}\]
then $\pi$ does not belong to $\Pi_u^1$. If the graph is
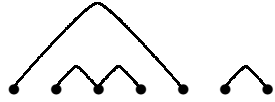
\[\begin{picture}(100,40)(0,0)
\courbe{0}{64}{0}{28}{32}{36}{32}
\courbe{16}{32}{0}{23}{24}{25}{8}
\courbe{32}{48}{0}{39}{40}{41}{8}
\courbe{80}{96}{0}{87}{88}{89}{8}
\put(-2,-3){$\bullet$}
\put(14,-3){$\bullet$}
\put(30,-3){$\bullet$}
\put(46,-3){$\bullet$}
\put(62,-3){$\bullet$}
\put(78,-3){$\bullet$}
\put(94,-3){$\bullet$}
\end{picture}\]
then $\pi$ lies in $\Pi_u^1$. Notice that these graphs generalize the notion of cup diagrams
(see \cite{Fung}, \cite{Melnikov-link-patterns}, \cite{Westbury}).

We characterize the components with a dense Jordan orbit as follows.

\begin{proposition}
\label{proposition-2}
The Jordan orbit ${\mathcal Z}_\pi\subset{\mathcal B}_u$ has maximal dimension
$\dim{\mathcal Z}_\pi=\dim{\mathcal B}_u$ if and only if $\pi\in\Pi_u^1$.
In particular, the mapping $\pi\mapsto {\mathcal K}^{T_\pi}$ is a one-to-one correspondence
between partitions $\pi\in\Pi_u^1$ and irreducible components of ${\mathcal B}_u$
which contain a dense Jordan orbit.
\end{proposition}

The proof is given in section \ref{section-2}.

Let $X,B_1,\ldots,B_m$ be algebraic varieties.
We recall from \cite{Fung} the notion of iterated fiber bundle, defined in the following inductive manner.
If $m=1$, then we say that $X$ is an iterated fiber bundle of base $B_1$ if there is an isomorphism
$X\stackrel{\mbox{\scriptsize $\sim$}}{\rightarrow} B_1$.
If $m>1$, we say that $X$ is an iterated fiber bundle of base $(B_1,\ldots,B_m)$
if there is a fiber bundle $X\rightarrow B_{m}$ whose fiber is an iterated fiber bundle of base
$(B_1,\ldots,B_{m-1})$.
For instance ${\mathcal B}^{(m)}$, the variety of complete flags of $\mathbb{C}^m$, is naturally an iterated fiber bundle
of base $(\mathbb{P}^{1},\ldots,\mathbb{P}^{m-1})$.
Then, it follows from formula (\ref{iso-Richardson}) that every Richardson component of ${\mathcal B}_u$
is an iterated fiber bundle of base $(\mathbb{P}^1,\ldots,\mathbb{P}^{\lambda^*_1-1},\ldots,
\mathbb{P}^1,\ldots,\mathbb{P}^{\lambda^*_s-1})$.
Or, in other words, if ${\mathcal K}^T\subset{\mathcal B}_u$ is a Bala-Carter component and $T^*$ denotes the
transposed tableau of $T$, then ${\mathcal K}^{T^*}$ is an iterated fiber bundle of base
$(\mathbb{P}^1,\ldots,\mathbb{P}^{\lambda_1-1},\ldots,
\mathbb{P}^1,\ldots,\mathbb{P}^{\lambda_r-1})$
(cf. section \ref{duality}).

More generally, we have the following result:

\begin{theorem}
\label{theorem-3}
Let $\lambda(u)=(\lambda_1\geq\ldots\geq\lambda_r)$ be the sizes of the Jordan blocks of $u$.
Let ${\mathcal K}^T\subset{\mathcal B}_u$ be an irreducible component,
associated to the tableau $T$.
Let ${\mathcal K}^{T^*}\subset{\mathcal B}_{u^*}$ be the component
corresponding to the transposed tableau $T^*$.
If ${\mathcal K}^T$ contains a dense Jordan orbit, then ${\mathcal K}^{T^*}$ is an iterated fiber bundle
of base $(\mathbb{P}^1,\ldots,\mathbb{P}^{\lambda_1-1},\ldots,
\mathbb{P}^1,\ldots,\mathbb{P}^{\lambda_r-1})$.
\end{theorem}

\begin{remark}
In the two-column case, every component contains a dense Jordan orbit (see \cite[\S 2]{Fresse-Melnikov-2}).
Theorem \ref{theorem-3} then implies that, in the two-row case, every component is an iterated fiber
bundle of base $(\mathbb{P}^1,\ldots,\mathbb{P}^1)$ ($\lambda_2$ terms).
In this manner, we retrieve a property which had been shown directly by F. Fung \cite{Fung}.
\end{remark}

\subsection{Outline}

The remainder of the paper comprises five parts.
Section \ref{section-2} is devoted to the proof of Proposition \ref{proposition-2}.
Using counting arguments, we provide an inductive estimation of the dimension of a Jordan orbit ${\mathcal Z}_\pi$
(formula (\ref{relation-codims})).
By induction, we derive that ${\mathcal Z}_\pi$ has the same dimension
as ${\mathcal B}_u$ if and only if $\pi\in\Pi_u^1$.
Proposition \ref{proposition-2} is obtained from this property.

The proof of Theorem \ref{theorem-3} is given in section \ref{section-6}.
It relies on the description of the standard tableaux associated to generalized Bala-Carter components
that Proposition \ref{proposition-2} provides.

Sections \ref{section-3}--\ref{section-5} are devoted to the proof of Theorem \ref{theorem-1}.
In section \ref{section-3}, some preliminary results are provided,
which are expressed as inductive criteria of singularity for Bala-Carter components
(Corollaries \ref{corollary-symmetry}-\ref{corollary-connected-component}).
In section \ref{section-4}, we prove the implication $(\Leftarrow)$ of Theorem \ref{theorem-1},
and to this end we construct two families of singular Bala-Carter components
and invoke the inductive arguments of the previous section.
In section \ref{section-5}, we prove the implication $(\Rightarrow)$ of Theorem \ref{theorem-1}. To do this,
again applying the preliminary results provided in section \ref{section-3},
we reduce the problem to check the smoothness of a single type of Bala-Carter components.
The proof that this component is indeed smooth (Proposition \ref{proposition-p-1-p-1})
is done by computing.

\subsubsection*{Notation.}
We set some conventional notation.
We denote by $\mathbf{S}_n$ the group of permutations of $\{1,\ldots,n\}$.
Let $\mathbb{P}^m$ be the projective space of dimension $m$ (i.e., the variety of linear lines
of $\mathbb{C}^{m+1}$).
If $V$ is a vector space, let $\mathrm{End}(V)$ be the space of endomorphisms of $V$
and let $GL(V)$ be the group of invertible endomorphisms.
If $A$ is a finite set, let $|A|$ denote its cardinal.
Moreover, if $A$ is composed of integers, let $\min A$ (resp. $\max A$) be its minimal (resp. maximal) element.
If $a,b$ are integers, let $[a,b]=\{a,\ldots,b\}$ be the integer interval between $a,b$,
and let $\,]a,b[\,=[a,b]\setminus\{a,b\}$.
If $X$ is an algebraic variety and $x\in X$, let $T_xX$ denote the tangent space of $X$
at the point $x$. If $Y\subset X$ is a subset, we denote by $\overline{Y}$ its closure
in the Zariski topology.
Other pieces of notation will be introduced in what follows.
The reader can find an index of the notation at the end of the article.

\section{Components with a dense Jordan orbit}

\label{section-2}

As in section \ref{section-background-statement}, we fix a nilpotent element $u\in\mathrm{End}(V)$.
We denote by $\lambda(u)=(\lambda_1\geq\ldots\geq\lambda_r)$ the sizes of the Jordan blocks of $u$,
by $Y(u)$ the Young diagram of rows of lengths $\lambda_1,\ldots,\lambda_r$,
and by $\lambda^*(u)=(\lambda^*_1\geq\ldots\geq\lambda^*_s)$ the lengths of the columns of $Y(u)$.
Let $Z_u=\{g\in GL(V):gug^{-1}=u\}$ be the stabilizer of $u$.
The purpose of this section is to show Proposition \ref{proposition-2}.

\subsection{Maximal dimensional Jordan orbits}

\label{section-3.1}

We consider an element $\pi\in\Pi_u$. Thus $\pi$ is a partition of $\{1,\ldots,n\}$
which can be written
$\pi=(I_1\sqcup \ldots\sqcup I_r)$, with $|I_j|=\lambda_j$ for all $j$.
Alternatively, $\pi$ can be seen as the map $\pi:\{1,\ldots,n\}\rightarrow \{\emptyset,1,\ldots,n\}$
with $\pi(i)=\emptyset$ if $i=\mathrm{min}(I_j)$, and $\pi(i)$ is the predecessor of $i$ in $I_j$
if $i\in I_j$, $i\not=\mathrm{min}(I_j)$.

We consider the $Z_u$-orbit ${\mathcal Z}_\pi\subset{\mathcal B}_u$. In this subsection, we determine under
which condition we have $\dim{\mathcal Z}_\pi=\dim{\mathcal B}_u$.
More precisely, our purpose is to show:

\begin{proposition}
\label{proposition-maximal-orbits}
We have $\dim {\mathcal Z}_\pi=\dim{\mathcal B}_u$ if and only if $\pi\in\Pi_u^1$.
\end{proposition}

We need some preliminary computations.
As a first step, we give the dimension of the group $Z_u$.

\begin{lemma}
\label{lemma-dim-centralizer}
We have $\dim Z_u=\displaystyle{\sum_{k=1}^s(\lambda_k^*)^2}$.
\end{lemma}

\medskip
\noindent {\em Proof.}
The group $Z_u$ is an open subset of the vector space
$\mathfrak{Z}_u=\{x\in\mathrm{End}(V):xu=ux\}$, hence $\dim Z_u=\dim \mathfrak{Z}_u$.
Let $W_1,\ldots,W_r$ be the Jordan blocks of $u$,
that is $V=W_1\oplus\ldots\oplus W_r$,
$\dim W_j=\lambda_j$ and $W_j=\langle e_j,u(e_j),\ldots,u^{\lambda_j-1}(e_j)\rangle$.
An element $x\in \mathfrak{Z}_u$ is
characterized by the images $x(e_j)$
for $j=1,\ldots,r$, and we necessarily have $x(e_j)\in\ker u^{\lambda_j}$.
Conversely, given $f_j\in \ker u^{\lambda_j}$ for every $j\in\{1,\ldots,r\}$,
there is a unique $x\in\mathfrak{Z}_u$ such that $x(e_j)=f_j$ for all $j$. It follows:
\[\dim\mathfrak{Z}_u=\sum_{j=1}^r\dim\ker u^{\lambda_j}=\sum_{k=1}^{s}|\{j:\lambda_j=k\}|.\dim\ker u^k.\]
Note that, for $k=1,\ldots,s$, $\dim\ker u^{k}=\lambda^*_1+\ldots+\lambda^*_k$.
In addition, $|\{j:\lambda_j=k\}|=\lambda_k^*-\lambda_{k+1}^*$
(by convention $\lambda^*_{s+1}=0$). It results that
\[\dim\mathfrak{Z}_u=\sum_{k=1}^{s}(\lambda_k^*-\lambda_{k+1}^*)(\lambda^*_1+\ldots+\lambda^*_k)=\sum_{k=1}^s(\lambda_k^*)^2.\]
This proves the lemma.
\hfill $\Box$

\medskip

Now, fix a $\pi$-basis $(e_1,\ldots,e_n)$
and let $F_0=(\langle e_1,\ldots,e_i\rangle)_{i=0,\ldots,n}\in{\mathcal Z}_\pi$ be the corresponding
adapted flag.
Then ${\mathcal Z}_\pi$ is the $Z_u$-orbit of $F_0$.
Let $Z_u^{F_0}=\{g\in Z_u:g(e_i)\in\langle e_1,\ldots,e_i\rangle\ \forall i=1,\ldots,n\}$ be the subgroup
of elements fixing $F_0$. Then
\begin{equation}
\label{dim-orbit}
\dim{\mathcal Z}_\pi=\dim Z_u-\dim Z_u^{F_0}.
\end{equation}
In the next step we determine $\dim Z_u^{F_0}$.
Set by convention $\pi(\emptyset)=\emptyset$ and $\emptyset<a$ for all $a\in\{1,\ldots,n\}$.
For $j=1,\ldots,r$, let $k_j=\max(I_j)$.
Define
\[A(\pi)=\{(i,j)\in \{1,\ldots,n\}\times\{1,\ldots,r\}:\pi^l(i)\leq \pi^l(k_j)\ \forall l\geq 0\}.\]
Then:

\begin{lemma}
\label{lemma-Api}
We have $\dim Z_u^{F_0}=|A(\pi)|$.
\end{lemma}

\noindent{\em Proof.}
The subgroup $Z_u^{F_0}$ is an open subset of the vector space
$\mathfrak{Z}_u^{F_0}=\{x\in\mathfrak{Z}_u:x(e_i)\in\langle e_1,\ldots,e_i\rangle\ \forall i=1,\ldots,n\}$,
hence $\dim Z_u^{F_0}=\dim \mathfrak{Z}_u^{F_0}$.

An element $x\in\mathfrak{Z}_u^{F_0}$ is determined by the images $x(e_{k_j})$
for $j=1,\ldots,r$. In addition, we must have
$x(e_{\pi^l(k_j)})\in\langle e_i:1\leq i\leq \pi^l(k_j)\rangle$ for all $j\in\{1,\ldots,r\}$, $l\geq 0$.
Recall that $x(e_{\pi^l(i)})=xu^l(e_i)=u^l(x(e_i))$ for all $i\in\{1,\ldots,n\}$, $l\geq 0$ (setting $e_\emptyset=0$).
Hence we must have $x(e_{k_j})\in\langle e_i:\pi^l(i)\leq \pi^l(k_j) \mbox{ for all }l\geq 0\rangle$.
Conversely, given $f_j\in\langle e_i:\pi^l(i)\leq \pi^l(k_j)\ \forall l\geq 0\rangle$ for $j=1,\ldots,r$,
there exists a unique $x\in\mathfrak{Z}_u^{F_0}$ such that $x(e_{k_j})=f_j$ for all $j$. It follows
\[\dim\mathfrak{Z}_u^{F_0}=\sum_{j=1}^r|\{i\in\{1,\ldots,n\}:\pi^l(i)\leq \pi^l(k_j)\ \forall l\geq 0\}|=|A(\pi)|.\]
The lemma is proved.
\hfill $\Box$

\medskip

Next, we establish an inductive estimate of the cardinal of the set $A(\pi)$.
We suppose $n\in I_{j_0}$ and we may suppose that $\lambda_{j_0}=|I_{j_0}|>|I_{j}|$ for $j>j_0$.
Let $I'_{j_0}=I_{j_0}\setminus\{n\}$, and let $I'_j=I_j$ for $j\not=j_0$.
Let $\lambda'_{j_0}=\lambda_{j_0}-1$ and let $\lambda'_j=\lambda_j$ for $j\not=j_0$.
Consider a nilpotent element $u'\in\mathrm{End}(\mathbb{C}^{n-1})$ of Jordan form $\lambda(u')=(\lambda'_1,\ldots,\lambda'_r)$.
Then $\pi'=(I'_1,\ldots,I'_r)$ belongs to $\Pi_{u'}$, it defines a $Z_{u'}$-orbit
${\mathcal Z}_{\pi'}\subset{\mathcal B}_{u'}$ and, in addition, it can be seen as a map
$\pi':\{1,\ldots,n-1\}\rightarrow\{\emptyset,1,\ldots,n-1\}$ which, in fact, is the restriction of $\pi$.

Let $Y(u')$ be the Young diagram associated to $\lambda(u')$, that is,
the sizes of its rows are $\lambda'_1,\ldots,\lambda'_r$.
Let $\lambda^*(u')=(\lambda'^*_1,\ldots,\lambda'^*_s)$ be the sizes of its columns.
Setting $l_0=\lambda_{j_0}$, we then have $\lambda'^*_{l_0}=j_0-1$ and
$\lambda'^*_l=\lambda^*_l$ for $l\not=l_0$. Therefore, using Lemma \ref{lemma-dim-centralizer}
and formula (\ref{relation-dimension-Bu}), we infer that
\begin{equation}
\label{inductive-dim-centralizer-Springer}
\dim Z_{u'}=\dim Z_u-2j_0+1\quad\mbox{and}\quad\dim {\mathcal B}_{u'}=\dim {\mathcal B}_{u}-j_0+1.
\end{equation}

We may consider the set $A(\pi')$ relative to $\pi'$.

\begin{lemma}
\label{lemma-Api-pi'} We have $|A(\pi)|\geq |A(\pi')|+j_0$, where
equality holds if and only if there are no $j>j_0$, $i\in I_j$ and
$l\geq 0$ such that $\pi^{l+1}(i)<\pi^{l+1}(n)<\pi^l(i)<\pi^l(n)$.
\end{lemma}

\noindent{\em Proof.}
For $j=1,\ldots,r$, we write $k'_j=\max(I'_j)$.
Notice that we possibly have $I'_{j_0}=\emptyset$.
In this case set by convention $k'_{j_0}=\emptyset$. The arguments remain valid in this case.
We always have $k'_{j_0}=\pi(n)$.

First, we show that $A(\pi')\subset A(\pi)$. Let $(i,j)\in A(\pi')$, that is,
$\pi'^l(i)\leq \pi'^l(k'_j)$ for all $l\geq 0$.
We have to show that $\pi^l(i)\leq \pi^l(k_j)$ for all $l\geq 0$.
If $j\not=j_0$, then $k'_j=k_j$ and, as $\pi'$ is the restriction of $\pi$ on $\{1,\ldots,n-1\}$,
the property holds and we infer that $(i,j)\in A(\pi)$.
Assume $j=j_0$.
Then $k'_j=\pi(n)$ and we have $\pi'^l(k'_j)=\pi^{l+1}(n)\leq \pi^{l}(n)$ for all $l\geq 0$.
Hence the property holds, and we get $(i,j)\in A(\pi)$.
Finally, we have obtained the desired inclusion $A(\pi')\subset A(\pi)$.

We have $(n,j_0)\in A(\pi)\setminus A(\pi')$.
Since $k_j<n$ for $j\not=j_0$,
by definition, there cannot be another element of the form $(n,j)$ in $A(\pi)$.

Let $j\in\{1,\ldots,j_0-1\}$. We have
$\pi^{\lambda_j-1}(k_j)>\pi^{\lambda_j-1}(k'_{j_0})=\emptyset$,
hence $(k_j,j_0)\notin A(\pi')$. Hence there is $i_j\in I_j$
minimal such that $(i_j,j_0)\notin A(\pi')$. By minimality, we
have $\pi^{l+1}(i_j)\leq\pi^l(k'_{j_0})=\pi^{l+1}(n)$ for all
$l\geq 0$, therefore, $(i_j,j_0)\in A(\pi)$. Thus, for $i\in I_j$,
we have $(i,j_0)\in A(\pi')\cap A(\pi)$ for $i<i_j$, $(i_j,j_0)\in
A(\pi)\setminus A(\pi')$, and necessarily $(i,j_0)\notin A(\pi)$
for $i>i_j$ (because otherwise, we would have $(i_j,j_0)\in
A(\pi')$).

At this stage, we have obtained that $|A(\pi)|\geq |A(\pi')|+j_0$,
with equality if and only if there are no $j>j_0$ and $i\in I_j$
such that $(i,j_0)\in A(\pi)\setminus A(\pi')$. It remains to show
that this condition is equivalent to the condition in the
statement of the lemma. Let us show this equivalence.

($\Leftarrow$) Suppose that there are $j>j_0$ and $i\in I_j$ such
that $(i,j_0)\in A(\pi)\setminus A(\pi')$. As $(i,j_0)\notin
A(\pi')$, there is $l\geq 0$ such that $\pi^l(i)>\pi^l(k'_{j_0})$,
that is $\pi^l(i)>\pi^{l+1}(n)$. As $(i,j_0)\in A(\pi)$, we have
$\pi^l(n)\geq \pi^l(i)$ and $\pi^{l+1}(n)\geq \pi^{l+1}(i)$. In
fact, $\pi^l(n)\not= \pi^l(i)$, because both are $\not=\emptyset$.
Moreover, $\pi^l(i)\not=\emptyset$ implies
$\pi^{l+1}(n)\not=\emptyset$, by the definition of $\pi$ and
because $|I_j|<|I_{j_0}|$, hence it also comes $\pi^{l+1}(n)\not=
\pi^{l+1}(i)$. Thus, we have found $j>j_0$, $i\in I_j$ and $l\geq
0$ such that $\pi^{l+1}(i)< \pi^{l+1}(n)<\pi^l(i)<\pi^l(n)$.

($\Rightarrow$) Suppose there are $j>j_0$, $i\in I_j$, $l\geq 0$
with $\pi^{l+1}(i)< \pi^{l+1}(n)<\pi^l(i)<\pi^l(n)$. In particular
$(i,j_0)\notin A(\pi')$. Then, there is $i_j\in I_j$ minimal such
that $(i_j,j_0)\notin A(\pi')$. As above, the minimality implies
$(i_j,j_0)\in A(\pi)$. Thus, we have found $j>j_0$ and $i_j\in
I_j$ such that $(i_j,j_0)\in A(\pi)\setminus A(\pi')$. \hfill
$\Box$

\medskip

Now we are ready to prove Proposition
\ref{proposition-maximal-orbits}.

\medskip

\noindent{\em Proof of Proposition
\ref{proposition-maximal-orbits}.} We reason by induction on
$n\geq 1$ with immediate initialization for $n=1$. Suppose that
the property holds for $n-1\geq 1$ and prove it for $n$. Combining
the relations (\ref{dim-orbit}) and
(\ref{inductive-dim-centralizer-Springer}) and Lemmas
\ref{lemma-Api}, \ref{lemma-Api-pi'}, we obtain
\begin{equation}
\label{relation-codims}
\dim{\mathcal B}_u-\dim{\mathcal Z}_\pi\geq \dim{\mathcal B}_{u'}-\dim{\mathcal Z}_{\pi'}
\end{equation}
with equality if and only if there are no $j>j_0$, $i\in I_j$ and
$l\geq 0$ such that $\pi^{l+1}(i)<\pi^{l+1}(n)<\pi^l(i)<\pi^l(n)$.

Let us interpret the last relation whenever it holds for some
$j>j_0$, $i\in I_j$, $l\geq 0$. If $\pi^{l+1}(i)>\emptyset$, then
it implies that $\pi$ has a crossing. If $\pi^{l+1}(i)=\emptyset$,
then $\pi^l(i)=\min(I_j)$ and we have $I_j<I_{j_0}$ though
$|I_j|<|I_{j_0}|$. In both cases, the relation implies
$\pi\notin\Pi_u^1$.

Suppose $\pi\in\Pi_u^1$. Then, the equality holds in relation
(\ref{relation-codims}). Let us show that $\pi'\in\Pi_{u'}^1$. It
is immediate that $\pi'$ has no crossing. Suppose $I'_j<I'_l$ with
$j,l\in\{1,\ldots,r\}$ and let us show $|I'_j|\geq |I'_l|$. If
$j_0\notin\{j,l\}$, then this follows from the fact that
$I'_j=I_j$, $I'_l=I_l$ and $\pi\in\Pi_u^1$. We cannot have
$I'_{j_0}<I'_l$, otherwise there is $i\in I'_l$ with
$\pi(i)<\pi(n)<i<n$ and $\pi$ has a crossing. If $I'_j<I'_{j_0}$,
then $I_j<I_{j_0}$ and we derive
$|I'_j|=|I_j|\geq|I_{j_0}|>|I'_{j_0}|$. Finally we have proved
$\pi'\in\Pi_{u'}^1$. By induction hypothesis, we have
$\dim{\mathcal Z}_{\pi'}=\dim{\mathcal B}_{u'}$. Therefore, we get
$\dim{\mathcal Z}_\pi=\dim{\mathcal B}_u$.

Conversely, suppose $\dim{\mathcal Z}_\pi=\dim{\mathcal B}_u$.
Then, necessarily, the equality holds in relation
(\ref{relation-codims}), and we have $\dim{\mathcal
Z}_{\pi'}=\dim{\mathcal B}_{u'}$. The latter fact implies
$\pi'\in\Pi_{u'}^1$ by induction hypothesis.
Because $\pi'$ has no crossing, a crossing of $\pi$, if it exists,
is of the form $\emptyset<\pi(i)<\pi(n)<i<n$. It follows from the
equality condition in (\ref{relation-codims}) that $i\in I_j$ with
$j<j_0$. There is $l\geq 1$ with $\pi^{l+1}(n)<\pi(i)<\pi^l(n)<i$.
In the case $\pi^{l+1}(n)\not=\emptyset$, we get that $\pi'$ has a
crossing. In the case $\pi^{l+1}(n)=\emptyset$, we get
$I'_{j_0}<I'_j$, though we have $|I'_{j_0}|<|I'_j|$. In both
cases, this contradicts the fact that $\pi'\in\Pi_{u'}^1$. Thus,
we obtain that $\pi$ has no crossing.
Next, we show: $I_j<I_l$ $\Rightarrow$ $|I_j|\geq|I_l|$. In the
case $j_0\notin\{j,l\}$, this follows from the fact that
$I'_j=I_j$, $I'_l=I_l$ and $\pi'\in\Pi_{u'}^1$. Note that we have
$I_{j_0}\not< I_j$ for all $j$. Thus, it remains to suppose
$I_j<I_{j_0}$ and then to show that $|I_j|\geq|I_{j_0}|$, that is
$j<j_0$. If $I'_j\not< I'_{j_0}$, then there is $i\in I_j$ such
that $\pi(i)<\pi(n)<i<n$. If $I'_j< I'_{j_0}$, then the fact that
$\pi'\in \Pi_{u'}^1$ implies $|I'_j|\geq |I'_{j_0}|$, hence there
are $i\in I_j$ and $l\geq 1$ such that
$\emptyset=\pi^{l+1}(i)<\pi^{l+1}(n)<\pi^l(i)=\min(I_j)<\pi^l(n)$.
In both cases, using the equality condition in
(\ref{relation-codims}), we infer that $j<j_0$.
Finally, we have shown $\pi\in\Pi_u^1$. The proof of the
proposition is then complete. \hfill $\Box$

\subsection{Proof of Proposition \ref{proposition-2}}

The first part of Proposition \ref{proposition-2} is provided by Proposition \ref{proposition-maximal-orbits}.
It remains to show the second part, that is:
the mapping $\pi\mapsto {\mathcal K}^{T_\pi}$ provides a one-to-one correspondence between elements
$\pi\in\Pi_u^1$ and components of ${\mathcal B}_u$ which contain a dense Jordan orbit.

First, suppose $\pi\in\Pi_u^1$.
From section \ref{Jordan-orbits}, we know that ${\mathcal Z}_\pi\subset{\mathcal K}^{T_\pi}$.
Then, by Proposition \ref{proposition-maximal-orbits}, we have that ${\mathcal Z}_\pi$ is indeed a dense
Jordan orbit of ${\mathcal K}^{T_\pi}$.
Hence, the mapping is well defined.

Conversely, suppose that the component ${\mathcal K}^T$ contains a dense Jordan orbit.
Recall that $\dim {\mathcal K}^T=\dim{\mathcal B}_u$, hence this orbit has the same dimension
as ${\mathcal B}_u$.
By Proposition \ref{proposition-maximal-orbits}, it is of the form ${\mathcal Z}_\pi$ with $\pi\in\Pi_u^1$.
As ${\mathcal Z}_\pi$ is dense in ${\mathcal K}^{T_\pi}$, we get ${\mathcal K}^T={\mathcal K}^{T_\pi}$.
Hence, the mapping is surjective.

Let $\pi,\pi'\in\Pi_u^1$ and suppose that ${\mathcal K}^{T_\pi}={\mathcal K}^{T_{\pi'}}$.
Thus, $\overline{{\mathcal Z}_\pi}=\overline{{\mathcal Z}_{\pi'}}$, which implies
${\mathcal Z}_\pi={\mathcal Z}_{\pi'}$ since two different orbits cannot have the same closure.
Let $\pi(j)=i\in\{\emptyset,1,\ldots,j-1\}$, and let $F=(V_0,...,V_n)\in{\mathcal Z}_\pi$.
Then $i$ is minimal such that $u(V_j)\subset V_i+u(V_{j-1})$.
Therefore, ${\mathcal Z}_\pi={\mathcal Z}_{\pi'}$ implies $\pi=\pi'$.
Hence, the mapping is injective. The proof of Proposition \ref{proposition-2}
is then complete.

\section{Inductive properties}

\label{section-3}

Let $\pi=(I_1,\ldots,I_r)\in\Pi_u$ and let ${\mathcal Z}_\pi$ be the corresponding Jordan orbit.
Following section \ref{Jordan-orbits}, the partition $\pi$ can be seen as a map
$\pi:\{1,\ldots,n\}\rightarrow\{\emptyset,1,\ldots,n\}$
such that for $i\in I_j$ we set $\pi(i)=\emptyset$ if $i=\min(I_j)$ and $\pi(i)$ to be the predecessor of $i$
in $I_j$ otherwise. Alternatively, $\pi$ is represented by a graph with $n$ vertices labeled by $1,\ldots,n$
displayed along a horizontal line, and with an arc between $i,j$ if $i=\pi(j)$.

We suppose in addition $\pi\in\Pi_u^1$. That is: $\pi$ has no crossing, and if $I_j\subset\,]\min(I_k),\max(I_k)[$
(i.e., $I_j<I_k$) then $|I_j|\geq |I_k|$.
Due to Proposition \ref{proposition-2}, we know that the orbit ${\mathcal Z}_\pi$ is dense in the
component ${\mathcal K}^{T_\pi}$.
In this section, our purpose is to establish some inductive singularity criteria for ${\mathcal K}^{T_\pi}$.
We consider $\pi'$ obtained from $\pi$ either by removing an extremal vertex ($1$ or $n$)
together with the corresponding arc, or by removing a whole connected component (from
the point of view of the graph representation). In both cases we show that the singularity of the component
${\mathcal K}^{T_{\pi'}}$ associated to $\pi'$ implies the singularity of ${\mathcal K}^{T_\pi}$.

We start by pointing out a symmetry property: if the graph of $\tilde\pi$ is the mirror of the graph of $\pi$,
then both components ${\mathcal K}^{T_\pi},{\mathcal K}^{T_{\tilde\pi}}$ are isomorphic.

\subsection{Symmetry}

For $j=1,\ldots,r$, we write $\tilde I_j=\{n-i+1:i\in I_j\}$. Then
$|\tilde I_j|=|I_j|$ for all $j$, and $\tilde\pi:=(\tilde I_1,\ldots,\tilde I_r)\in\Pi_u$
is another partition of $\{1,\ldots,n\}$. The graph representation of $\tilde\pi$
is symmetric to the graph representation of $\pi$ (that is, it is its mirror reflection).
Example:
\[\pi=\{1,5,6\}\sqcup\{2,3,4\}\sqcup\{8,9\}\sqcup\{7\}
\ \rightarrow \
\tilde\pi=\{4,5,9\}\sqcup\{6,7,8\}\sqcup\{1,2\}\sqcup\{3\}
\]
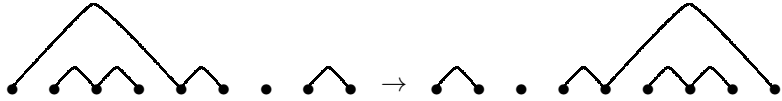
\[
\begin{picture}(130,45)(0,-5)
\courbe{0}{64}{0}{29}{31}{35}{32}
\courbe{16}{32}{0}{23}{24}{25}{8}
\courbe{32}{48}{0}{39}{40}{41}{8}
\courbe{64}{80}{0}{71}{72}{73}{8}
\courbe{112}{128}{0}{119}{120}{121}{8}
\put(-2,-3){$\bullet$}
\put(14,-3){$\bullet$}
\put(30,-3){$\bullet$}
\put(46,-3){$\bullet$}
\put(62,-3){$\bullet$}
\put(78,-3){$\bullet$}
\put(94,-3){$\bullet$}
\put(110,-3){$\bullet$}
\put(126,-3){$\bullet$}
\end{picture}
\quad { }^{\displaystyle\rightarrow} \quad
\begin{picture}(130,45)(0,-5)
\courbe{0}{16}{0}{7}{8}{9}{8}
\courbe{48}{64}{0}{55}{56}{57}{8}
\courbe{64}{128}{0}{93}{96}{99}{32}
\courbe{80}{96}{0}{87}{88}{89}{8}
\courbe{96}{112}{0}{103}{104}{105}{8}
\put(-2,-3){$\bullet$}
\put(14,-3){$\bullet$}
\put(30,-3){$\bullet$}
\put(46,-3){$\bullet$}
\put(62,-3){$\bullet$}
\put(78,-3){$\bullet$}
\put(94,-3){$\bullet$}
\put(110,-3){$\bullet$}
\put(126,-3){$\bullet$}
\end{picture}
\]
As we suppose $\pi\in\Pi_u^1$, it is clear that $\tilde\pi\in\Pi_u^1$.
Thus, the partition $\tilde\pi$ defines a component ${\mathcal K}^{T_{\tilde\pi}}\subset{\mathcal B}_u$
which is the closure of the Jordan orbit ${\mathcal Z}_{\tilde\pi}$. We show:

\begin{proposition}
\label{proposition-symmetry}
The components ${\mathcal K}^{T_\pi}$ and ${\mathcal K}^{T_{\tilde\pi}}$ are isomorphic.
\end{proposition}

\noindent{\em Proof.} Let $\tilde V$ be the vector space of linear forms $\phi:V\rightarrow \mathbb{C}$,
and let $\tilde u:\tilde V\rightarrow \tilde V$, $\phi\mapsto\phi\circ\tilde u$ be the
dual map of $u$. Then $\tilde u$ is nilpotent and has the same Jordan form as $u$.
Let ${\mathcal B}_{\tilde u}$ be the variety of $\tilde u$-stable flags of $\tilde V$.
If $W\subset V$ is a subspace, we set $W^\perp=\{\phi\in\tilde V:\phi(w)=0\ \forall w\in W\}$.
Then, $W$ is $u$-stable if and only if $W^\perp$ is $\tilde u$-stable,
and the map $\Phi:{\mathcal B}_u\rightarrow {\mathcal B}_{\tilde u}$,
$(V_0,\ldots,V_n)\mapsto (V_n^\perp,\ldots,V_0^\perp)$ is an isomorphism of algebraic varieties.
In particular, $\Phi({\mathcal K}^{T_\pi})\subset{\mathcal B}_{\tilde u}$ is an irreducible component.

Let $Z_{\tilde u}=\{g\in GL(\tilde V):g\tilde ug^{-1}=\tilde u\}$ be the stabilizer of $\tilde u$.
We denote by $\tilde{\mathcal Z}_{\tilde\pi}\subset{\mathcal B}_{\tilde u}$ the Jordan $Z_{\tilde u}$-orbit
associated to $\tilde\pi$, and by $\tilde{\mathcal K}^{T_{\tilde \pi}}\subset{\mathcal B}_{\tilde u}$
the irreducible component associated to the tableau $T_{\tilde \pi}$, which is actually the closure of
$\tilde{\mathcal Z}_{\tilde\pi}$. It is easy to see that $\Phi({\mathcal Z}_\pi)=\tilde{\mathcal Z}_{\tilde\pi}$.
Thus, $\Phi({\mathcal K}^{T_\pi})=\tilde{\mathcal K}_{\tilde\pi}$
so that $\Phi$ restricts to an isomorphism from ${\mathcal K}^{T_\pi}$ to $\tilde{\mathcal K}^{T_{\tilde\pi}}$.

We distinguish between $\tilde{\mathcal K}^{T_{\tilde\pi}}\subset{\mathcal B}_{\tilde u}$,
the component associated to the standard tableau $T_{\tilde\pi}$ in ${\mathcal B}_{\tilde u}$,
and ${\mathcal K}^{T_{\tilde\pi}}$, the component associated to the same tableau, but in ${\mathcal B}_u$.
A component only depends up to isomorphism on the standard tableau it is associated to,
hence these two components are actually isomorphic.
We deduce that ${\mathcal K}^{T_\pi}$ and ${\mathcal K}^{T_{\tilde\pi}}$ are isomorphic.
\hfill $\Box$

\begin{remark}
We can see that the tableau $T_{\tilde \pi}$ is in fact the image of the tableau $T_\pi$ by the classical
Sch$\ddot{\mbox{u}}$tzenberger involution.
\end{remark}

\subsection{Removing an extremal point}

\label{section-extremal-point}

There is $j_0\in\{1,\ldots,r\}$ such that $n\in I_{j_0}$.
We may assume that
$\lambda_{j_0}=|I_{j_0}|>|I_{j}|$ for all $j>j_0$.
As in section \ref{section-3.1},
let $\lambda'_{j_0}=\lambda_{j_0}-1$ and let $\lambda'_j=\lambda_j$ for $j\not=j_0$,
and consider a nilpotent element $u'\in\mathrm{End}(\mathbb{C}^{n-1})$ of Jordan form $\lambda(u')=(\lambda'_1,\ldots,\lambda'_r)$.
Let $I'_{j_0}=I_{j_0}\setminus\{n\}$, and let $I'_j=I_j$ for $j\not=j_0$, then
$\pi':=(I'_1,\ldots,I'_r)$ belongs to $\Pi_{u'}$ and it defines a $Z_{u'}$-orbit
${\mathcal Z}_{\pi'}\subset{\mathcal B}_{u'}$.
At the level of the graph representation, the graph of $\pi'$ is obtained from the graph of $\pi$
by removing the last vertex together with the possible upcoming arc, for instance:
\[\pi=\{1,5,6\}\sqcup\{2,3,4\}\sqcup\{8,9\}\sqcup\{7\}
\ \rightarrow \
\pi'=\{1,5,6\}\sqcup\{2,3,4\}\sqcup\{8\}\sqcup\{7\}\ \,
\]
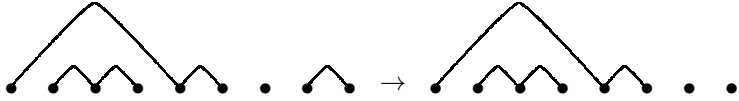
\[
\begin{picture}(130,45)(0,-5)
\courbe{0}{64}{0}{29}{32}{35}{32}
\courbe{16}{32}{0}{23}{24}{25}{8}
\courbe{32}{48}{0}{39}{40}{41}{8}
\courbe{64}{80}{0}{71}{72}{73}{8}
\courbe{112}{128}{0}{119}{120}{121}{8}
\put(-2,-3){$\bullet$}
\put(14,-3){$\bullet$}
\put(30,-3){$\bullet$}
\put(46,-3){$\bullet$}
\put(62,-3){$\bullet$}
\put(78,-3){$\bullet$}
\put(94,-3){$\bullet$}
\put(110,-3){$\bullet$}
\put(126,-3){$\bullet$}
\end{picture}
\quad { }^{\displaystyle\rightarrow} \quad
\begin{picture}(130,45)(0,-5)
\courbe{0}{64}{0}{29}{32}{35}{32}
\courbe{16}{32}{0}{23}{24}{25}{8}
\courbe{32}{48}{0}{39}{40}{41}{8}
\courbe{64}{80}{0}{71}{72}{73}{8}
\put(-2,-3){$\bullet$}
\put(14,-3){$\bullet$}
\put(30,-3){$\bullet$}
\put(46,-3){$\bullet$}
\put(62,-3){$\bullet$}
\put(78,-3){$\bullet$}
\put(94,-3){$\bullet$}
\put(110,-3){$\bullet$}
\end{picture}
\]
It is clear that $\pi'\in\Pi_{u'}^1$, hence the orbit ${\mathcal Z}_{\pi'}$ is
dense in the component ${\mathcal K}^{T_{\pi'}}\subset{\mathcal B}_{u'}$.
We have:

\begin{proposition}
\label{proposition-removing-n}
If the component ${\mathcal K}^{T_{\pi'}}$ is singular, then ${\mathcal K}^{T_\pi}$ is singular.
Suppose moreover that $n\in I_{j_0}$ with $|I_{j_0}|=\max\{|I_j|:1\leq j\leq r\}$,
then, ${\mathcal K}^{T_{\pi'}}$ is singular if and only if ${\mathcal K}^{T_\pi}$ is singular.
\end{proposition}

\noindent
{\em Proof.}
Notice that it follows from the definition of the tableaux $T_\pi$, $T_{\pi'}$ in section \ref{Jordan-orbits} that
$T_{\pi'}$ is obtained from $T_\pi$ simply by deleting the box $n$.
Moreover, under the assumption $|I_{j_0}|=\max\{|I_j|:1\leq j\leq r\}$,
$n$ lies in the last column of $T_\pi$.
Then the result follows from \cite[Theorem 2.1]{Fresse-Melnikov-1}, which says more generally
that, if $T$ is a standard tableau and $T'$ is the subtableau obtained by deleting $n$,
then the singularity of the component ${\mathcal K}^{T'}$ implies the singularity of ${\mathcal K}^T$,
and moreover, if $n$ lies in the last column of $T$,
then the singularity of ${\mathcal K}^{T'}$ is equivalent to the singularity of ${\mathcal K}^T$.
\hfill $\Box$

\medskip

We have just studied the situation where we remove from $\pi$ the last number $n$,
the situation where we remove the first number $1$ is quite similar.
Let $j_1\in\{1,\ldots,r\}$ be such that $1\in I_{j_1}$
and we may suppose $\lambda_{j_1}=|I_{j_1}|>|I_{j}|$ for $j>j_1$.
Let $\lambda''_{j_1}=\lambda_{j_1}-1$ and $\lambda''_j=\lambda_j$ for $j\not=j_1$.
Fix $u''\in\mathrm{End}(\mathbb{C}^{n-1})$ nilpotent of Jordan form
$\lambda(u'')=(\lambda''_1,\ldots,\lambda''_r)$.
Let $I''_{j_1}=\{i-1:i\in I_{j_1}\setminus\{1\}\}$ and
$I''_j=\{i-1:i\in I_j\}$ for $j\not=j_1$,
then $\pi''=(I''_1,\ldots,I''_r)$ belongs to $\Pi_{u''}$
and it defines a $Z_{u''}$-orbit ${\mathcal Z}_{\pi''}\subset{\mathcal B}_{u''}$.
At the level of graphs, the graph of $\pi''$ is obtained from the graph of $\pi$ by
removing the first vertex and the possible upcoming arc.
\[\pi=\{1,5,6\}\sqcup\{2,3,4\}\sqcup\{8,9\}\sqcup\{7\}
\ \rightarrow \
\pi''=\{1,2,3\}\sqcup\{4,5\}\sqcup\{7,8\}\sqcup\{6\}\ \,
\]
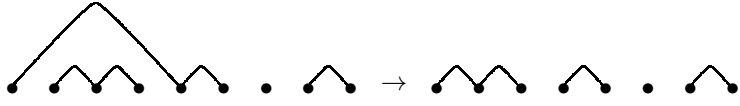
\[
\begin{picture}(130,45)(0,-5)
\courbe{0}{64}{0}{29}{32}{35}{32}
\courbe{16}{32}{0}{23}{24}{25}{8}
\courbe{32}{48}{0}{39}{40}{41}{8}
\courbe{64}{80}{0}{71}{72}{73}{8}
\courbe{112}{128}{0}{119}{120}{121}{8}
\put(-2,-3){$\bullet$}
\put(14,-3){$\bullet$}
\put(30,-3){$\bullet$}
\put(46,-3){$\bullet$}
\put(62,-3){$\bullet$}
\put(78,-3){$\bullet$}
\put(94,-3){$\bullet$}
\put(110,-3){$\bullet$}
\put(126,-3){$\bullet$}
\end{picture}
\quad { }^{\displaystyle\rightarrow} \quad
\begin{picture}(130,45)(16,-5)
\courbe{16}{32}{0}{23}{24}{25}{8}
\courbe{32}{48}{0}{39}{40}{41}{8}
\courbe{64}{80}{0}{71}{72}{73}{8}
\courbe{112}{128}{0}{119}{120}{121}{8}
\put(14,-3){$\bullet$}
\put(30,-3){$\bullet$}
\put(46,-3){$\bullet$}
\put(62,-3){$\bullet$}
\put(78,-3){$\bullet$}
\put(94,-3){$\bullet$}
\put(110,-3){$\bullet$}
\put(126,-3){$\bullet$}
\end{picture}
\]
Clearly $\pi''\in\Pi_{u''}^1$, hence ${\mathcal Z}_{\pi''}$
is dense in the component ${\mathcal K}^{T_{\pi''}}\subset{\mathcal B}_{u''}$.
Combining Propositions \ref{proposition-symmetry} and \ref{proposition-removing-n},
we derive:

\begin{proposition}
If the component ${\mathcal K}^{T_{\pi''}}$ is singular, then ${\mathcal K}^{T_\pi}$ is singular.
Suppose moreover $1\in I_{j_1}$ with $|I_{j_1}|=\max\{|I_j|:1\leq j\leq r\}$,
then, ${\mathcal K}^{T_{\pi''}}$ is singular if and only if ${\mathcal K}^{T_\pi}$ is singular.
\end{proposition}

\subsection{Removing a connected component}

We consider the situation of a partition obtained from $\pi$ by removing a term in the partition.
We suppose $\pi=(I_1,\ldots,I_r)\in\Pi_u$ with $|I_j|=\lambda_j$ for all $j$.
Fix $\hat{\jmath}\in\{1,\ldots,r\}$ and fix a nilpotent endomorphism $\hat{u}\in\mathrm{End}(\mathbb{C}^{n-\lambda_{\hat\jmath}})$
of Jordan form $\lambda(\hat{u})=(\lambda_1,\ldots,\lambda_{\hat{\jmath}-1},
\lambda_{\hat{\jmath}+1},\ldots,\lambda_r)$.
Write $\{1,\ldots,n\}\setminus I_{\hat\jmath}=\{k_1<\ldots<k_{n-\lambda_{\hat\jmath}}\}$. 
For $j\not=\hat\jmath$, let $\hat{I}_j=\{i:k_i\in I_j\}$. Then
$\hat\pi=(\hat{I}_1,\ldots,\hat{I}_{\hat\jmath-1},\hat{I}_{\hat\jmath+1},\ldots,\hat{I}_r)$ is an element of $\Pi_{\hat u}$.
At the level of graphs, the graph of $\hat\pi$ is obtained from the graph of $\pi$ by removing a connected
component, the one which corresponds to the subset $I_{\hat\jmath}$. For example:
\[\pi=\{1,5,6\}\sqcup\underbrace{\{2,3,4\}}_{\mbox{\scriptsize to remove}}\sqcup\{8,9\}\sqcup\{7\}
\ \rightarrow \
\hat\pi=\{1,2,3\}\sqcup\{5,6\}\sqcup\{4\}\ \ \ \ \ \ \
\]
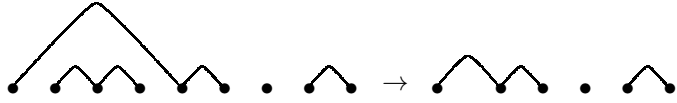
\[
\begin{picture}(130,36)(0,-5)
\courbe{0}{64}{0}{29}{32}{35}{32}
\courbe{16}{32}{0}{23}{24}{25}{8}
\courbe{32}{48}{0}{39}{40}{41}{8}
\courbe{64}{80}{0}{71}{72}{73}{8}
\courbe{112}{128}{0}{119}{120}{121}{8}
\put(-2,-3){$\bullet$}
\put(14,-3){$\bullet$}
\put(30,-3){$\bullet$}
\put(46,-3){$\bullet$}
\put(62,-3){$\bullet$}
\put(78,-3){$\bullet$}
\put(94,-3){$\bullet$}
\put(110,-3){$\bullet$}
\put(126,-3){$\bullet$}
\end{picture}
\quad { }^{\displaystyle\rightarrow} \quad
\begin{picture}(115,36)(8,-5)
\courbe{8}{32}{0}{18}{20}{22}{12}
\courbe{32}{48}{0}{39}{40}{41}{8}
\courbe{80}{96}{0}{87}{88}{89}{8}
\put(6,-3){$\bullet$}
\put(30,-3){$\bullet$}
\put(46,-3){$\bullet$}
\put(62,-3){$\bullet$}
\put(78,-3){$\bullet$}
\put(94,-3){$\bullet$}
\end{picture}
\]
As we suppose $\pi\in\Pi_u^1$, it is clear that $\hat\pi\in\Pi_{\hat u}^1$. Thus,
the partition $\hat\pi$ defines a component ${\mathcal K}^{T_{\hat\pi}}\subset{\mathcal B}_{\hat u}$
which is the closure of the Jordan orbit ${\mathcal Z}_{\hat\pi}$. We show:

\begin{proposition}
\label{proposition-remove-connected-component}
If the component ${\mathcal K}^{T_{\hat\pi}}$ is singular, then ${\mathcal K}^{T_\pi}$ is singular.
\end{proposition}

\noindent
{\em Proof.} Suppose that ${\mathcal K}^{T_\pi}$ is smooth and let us show
that ${\mathcal K}^{T_{\hat\pi}}$ is smooth.
Recall that $\pi$ can be seen as a map $\pi:\{1,\ldots,n\}\rightarrow\{\emptyset,1,\ldots,n\}$.
We fix a $\pi$-basis $(e_1,\ldots,e_n)$, that is $u(e_i)=e_{\pi(i)}$ for all $i$,
with by convention $e_\emptyset=0$.
Then $W_j=\langle e_i:i\in I_j\rangle$ (for $j=1,\ldots,r$) are the Jordan blocks of $u$.
For $t\in\mathbb{C}^*$, we consider $h_t:V\rightarrow V$
defined by $(h_t)_{|W_j}=\mathrm{id}_{W_j}$ for $j\not=\hat\jmath$ and
$(h_t)_{|W_{\hat\jmath}}=t.\mathrm{id}_{W_{\hat\jmath}}$.
Then, $H=\{h_t:t\in\mathbb{C}^*\}$ is a subtorus of rank one of the group $Z_u$.
Therefore, it acts on ${\mathcal B}_u$ and stabilizes its irreducible components,
in particular it stabilizes ${\mathcal K}^{T_\pi}$.
As we suppose that ${\mathcal K}^{T_\pi}$ is smooth, 
we know that the fixed point set
$({\mathcal K}^{T_\pi})^H=\{F\in {\mathcal K}^{T_\pi}: h_tF=F\ \forall t\in\mathbb{C}^*\}$ is a smooth subvariety
(cf., \cite{Bialynicki-Birula}).
To prove the proposition, it is then sufficient to show that ${\mathcal K}^{T_{\hat\pi}}$
is isomorphic to a connected component of
$({\mathcal K}^{T_\pi})^H$. This is what we do in the following.

Let $\hat{V}=\bigoplus_{j\not=\hat\jmath}W_j$ and write $\hat u$ the restriction of $u$ to $\hat V$
(this accords with the previous notation $\hat u$). Let ${\mathcal B}_{\hat u}$ be the variety of $\hat u$-stable
complete flags of $\hat V$. We construct a map $\Phi:{\mathcal B}_{\hat u}\rightarrow {\mathcal B}_u$ as follows.
For $i=0,\ldots,n$, set $a_i=|\{1,\ldots,i\}\cap I_{\hat\jmath}|$.
Write $I_{\hat\jmath}=\{\hat k_1<\ldots<\hat k_{\lambda_{\hat\jmath}}\}$,
and for $a\in\{0,\ldots,\lambda_{\hat\jmath}\}$ set $W_{\hat\jmath,a}=\langle e_{\hat k_1},\ldots,e_{\hat k_a}\rangle$.
Hence, $(W_{\hat\jmath,0}\subset W_{\hat\jmath,1}\subset\ldots\subset W_{\hat\jmath,\lambda_{\hat\jmath}})$
is the unique $u$-stable complete flag of $W_{\hat\jmath}$.
Define $\Phi$ by
\[\Phi:{\mathcal B}_{\hat u}\rightarrow {\mathcal B}_u,\ (\hat V_0,\ldots,\hat V_{n-\lambda_{\hat\jmath}})
\mapsto(\hat V_{i-a_i}\oplus W_{\hat\jmath,a_i}:i=0,\ldots,n).\]
The map $\Phi$ is well defined and algebraic, and its image lies in the fixed point set
$({\mathcal B}_u)^H$ for the action of $H$ on ${\mathcal B}_u$.
Notice that a flag $F=(V_0,\ldots,V_n)\in ({\mathcal B}_u)^H$ satisfies
$V_i=(V_i\cap\hat{V})\oplus (V_i\cap W_{\hat\jmath})$ for all $i$.
Then, the image of $\Phi$ can be characterized as follows:
\begin{eqnarray*}
\Phi({\mathcal B}_{\hat u}) &\!\!\!=\!\!\!& \{F=(V_0,\ldots,V_n)\in ({\mathcal B}_u)^H:\dim V_i\cap W_{\hat\jmath}=a_i \ \forall i=0,\ldots,n\} \\
 &\!\!\!=\!\!\!& \{F\in ({\mathcal B}_u)^H:\dim V_i\cap W_{\hat\jmath}\leq a_i\mbox{ and }
 \dim V_i\cap \hat V\leq i-a_i\ \forall i=0,\ldots,n\} \\
 &\!\!\!=\!\!\!& \{F\in ({\mathcal B}_u)^H:\dim V_i\cap W_{\hat\jmath}\geq a_i\mbox{ and }
 \dim V_i\cap \hat V\geq i-a_i\ \forall i=0,\ldots,n\}.
\end{eqnarray*}
Hence $\Phi({\mathcal B}_{\hat u})$ is open and closed in $({\mathcal B}_u)^H$,
and the map $\Phi$ is an open and closed immersion.

The intersection $\Phi({\mathcal B}_{\hat u})\cap{\mathcal K}^{T_\pi}$ is also open and closed in
$({\mathcal K}^{T_\pi})^H$, hence it is a union
$C_1\sqcup\ldots\sqcup C_m$ of connected
components of $({\mathcal K}^{T_\pi})^H$.
It is easy to see that, by construction, we have the inclusion
$\Phi({\mathcal Z}_{\hat \pi})\subset {\mathcal Z}_\pi$,
which implies $\Phi({\mathcal K}^{T_{\hat\pi}})\subset ({\mathcal K}^{T_\pi})^H$.
As $\Phi({\mathcal K}^{T_{\hat\pi}})$ is irreducible, we have $\Phi({\mathcal K}^{T_{\hat\pi}})\subset C_l$
for some $l\in\{1,\ldots,m\}$. Since $\Phi({\mathcal K}^{T_{\hat\pi}})$ has the same dimension as $\Phi({\mathcal B}_{\hat u})$
and since $C_l$ is irreducible (because $({\mathcal K}^{T_\pi})^H$ is smooth),
we actually have $\Phi({\mathcal K}^{T_{\hat\pi}})= C_l$.
Therefore, $\Phi$ restricts to an isomorphism between the component ${\mathcal K}^{T_{\hat\pi}}$
and a connected component of $({\mathcal K}^{T_\pi})^H$. The proof of the proposition
is then complete.
\hfill $\Box$

\subsection{At the level of Bala-Carter components}

In the previous subsections, we have provided inductive criteria
for the singularity of
an irreducible component of the form ${\mathcal K}^{T_\pi}\subset {\mathcal B}_u$
with $\pi\in\Pi_u^1$, that is, more general than a component of Bala-Carter type.
In this subsection, we translate the previous criteria for a Bala-Carter component.

We consider an element $\pi\in \Lambda_u$,
that is, a sequence $\pi=(\pi_1,\ldots,\pi_r)$ such that the numbers $\pi_1,\ldots,\pi_r$
coincide up to ordering with the Jordan block sizes $\lambda_1,\ldots,\lambda_r$.
This sequence defines the Bala-Carter component ${\mathcal K}_\pi^{\mathrm{BC}}$.
According to subsection \ref{Jordan-orbits}, $\pi$ can be seen as an element of the set $\Pi_u^0$,
that we also denote by $\pi$,
and the component ${\mathcal K}_\pi^{\mathrm{BC}}$ then coincides with the component ${\mathcal K}^{T_\pi}$.
Through this identification, we derive the following corollaries from Propositions
\ref{proposition-symmetry}--\ref{proposition-remove-connected-component}.

Nota: a sequence of nonnegative integers defines a Bala-Carter component in the appropriate Springer fiber.
In the following statements, we do not precise the underlying Springer fibers in which the components
are imbedded.

\begin{corollary}
\label{corollary-symmetry}
Let $\tilde\pi=(\pi_r,\pi_{r-1},\ldots,\pi_1)$.
Then the Bala-Carter components ${\mathcal K}^{\mathrm{BC}}_\pi$ and ${\mathcal K}^{\mathrm{BC}}_{\tilde\pi}$
are isomorphic.
\end{corollary}

\begin{corollary}
\label{corollary-extremal-points}
Let $\pi'=(\pi_1,\ldots,\pi_{r-1},\pi_r-1)$.
If the component ${\mathcal K}^{\mathrm{BC}}_{\pi'}$ is singular, then ${\mathcal K}^{\mathrm{BC}}_\pi$ is singular.
Suppose moreover $\pi_r=\max\{\pi_j:j=1,\ldots,r\}$, then
${\mathcal K}^{\mathrm{BC}}_{\pi'}$ is singular if and only if ${\mathcal K}^{\mathrm{BC}}_\pi$ is singular.

Let $\pi''=(\pi_1-1,\pi_2,\ldots,\pi_r)$.
Likewise, if ${\mathcal K}^{\mathrm{BC}}_{\pi''}$ is singular, then ${\mathcal K}^{\mathrm{BC}}_\pi$ is singular.
Moreover, if $\pi_1=\max\{\pi_j:j=1,\ldots,r\}$, then
${\mathcal K}^{\mathrm{BC}}_{\pi''}$ is singular if and only if ${\mathcal K}^{\mathrm{BC}}_\pi$ is singular.
\end{corollary}

\begin{corollary}
\label{corollary-connected-component}
Let $\hat\jmath\in\{1,\ldots,r\}$ and $\hat\pi=(\pi_1,\ldots,\pi_{\hat\jmath-1},\pi_{\hat\jmath+1},\ldots,\pi_r)$.
If the component ${\mathcal K}^{\mathrm{BC}}_{\hat\pi}$ is singular, then ${\mathcal K}^{\mathrm{BC}}_\pi$ is singular.
\end{corollary}

\section{Singular Bala-Carter components}

\label{section-4}

Let $\pi\in\Lambda_u$, that is, a sequence $\pi=(\pi_1,\ldots,\pi_r)$ which coincides
up to ordering
with the sequence $\lambda(u)=(\lambda_1\geq\ldots\geq\lambda_r)$ of the Jordan block sizes of $u$.
We associate to $\pi$ the Bala-Carter component ${\mathcal K}_\pi^{\mathrm{BC}}\subset{\mathcal B}_u$.
The purpose of this section is to show the implication $(\Leftarrow)$
of Theorem \ref{theorem-1}, namely:

\begin{proposition}
\label{proposition-theorem-1-if}
If $\pi\geq(1,2,2,1)$ or $\pi\geq(2,3,2)$, then ${\mathcal K}_\pi^{\mathrm{BC}}$ is singular.
\end{proposition}

First, we construct two families of singular Bala-Carter components.
Then, we derive the proposition by using the results in the previous section.

\subsection{Singular components of the form ${\mathcal K}_\pi^{\mathrm{BC}}$ with $\pi=(1,p,q,1)$}

In \cite{Vargas}, J.A. Vargas gives an example of a singular component, which is
${\mathcal K}_\pi^{\mathrm{BC}}$ for $\pi=(1,2,2,1)$.
Here, we provide a family of singular components which generalizes this example:

\begin{lemma}
\label{lemma-singularity-1-p-q-1}
Let $p\geq 2$, $q\geq 2$, and $\pi=(1,p,q,1)$. Then the component ${\mathcal K}_\pi^{\mathrm{BC}}$ is singular.
\end{lemma}

\noindent {\em Proof.}
Due to Corollary \ref{corollary-symmetry}, we may assume $p\geq q$.
Moreover we may assume $p>2$, otherwise $\pi=(1,2,2,1)$ and we already know that the component
${\mathcal K}_\pi^{\mathrm{BC}}$ is singular by Vargas's result.
The component lies in ${\mathcal B}_u$
where $u\in\mathrm{End}(\mathbb{C}^{p+q+2})$ has Jordan form $\lambda(u)=(p,q,1,1)$.
Due to formula (\ref{relation-dimension-Bu}), we have $\dim {\mathcal K}_\pi^{\mathrm{BC}}=\dim{\mathcal B}_u=q+5$.

We fix a Jordan basis $(e_1,\ldots,e_{p+q+2})$ such that $u$ acts on the basis according to the following
picture
\[
\begin{array}{ccccccccccccccccccccccccccccccccc}
0 & \leftarrow & e_1 & \leftarrow & e_3 & \leftarrow & \cdots & \leftarrow & e_{p} & \leftarrow & e_{p+q+1} \\
0 & \leftarrow & e_2 & \leftarrow & e_{p+3} & \leftarrow & \cdots & e_{p+q} & \leftarrow & e_{p+q+2} \\
0 & \leftarrow & e_{p+1} \\
0 & \leftarrow & e_{p+2}
\end{array}
\]
Note that, if $q=2$, then the sequence $(e_{p+3},\ldots,e_{p+q})$ disappears from the second line,
and the picture becomes
\[
\begin{array}{ccccccccccccccccccccccccccccccccc}
0 & \leftarrow & e_1 & \leftarrow & e_3 & \leftarrow & \cdots & \leftarrow & e_{p} & \leftarrow & e_{p+3} \\
0 & \leftarrow & e_2 & \leftarrow & e_{p+4} \\
0 & \leftarrow & e_{p+1} \\
0 & \leftarrow & e_{p+2}
\end{array}
\]
We consider the flag $F_0=(\langle e_1,\ldots,e_i\rangle)_{i=0,\ldots,p+q+2}$ adapted to the basis,
and the purpose is to show that $F_0$ is a singular point of the component ${\mathcal K}_\pi^{\mathrm{BC}}$.
To do this, we first check that $F_0$ lies in ${\mathcal K}_\pi^{\mathrm{BC}}$ and then we construct
$q+6$ linearly independent vectors of the tangent space of ${\mathcal K}_\pi^{\mathrm{BC}}$ at the point $F_0$.

Consider $B=\{g\in GL(V):ge_i\in\langle e_i,\ldots,e_{p+q+2}\rangle\}$
the subgroup of lower triangular matrices in the basis.
Then the orbit $\Omega=B F_0$ is an affine open neighborhood of $F_0$
in the variety of complete flags of $\mathbb{C}^{p+q+2}$.
Consider $U\subset B$ the subgroup of unipotent matrices
and $\mathfrak{n}=\{g-I:g\in U\}$ its Lie algebra.
The map $\varphi:\mathfrak{n}\rightarrow \Omega$, $x\mapsto (I+x)F_0$ is an isomorphism of affine varieties
and we consider the subvariety $\varphi^{-1}(\Omega\cap{\mathcal K}_\pi^{\mathrm{BC}})\subset\mathfrak{n}$.
For $1\leq k<l\leq p+q+2$, let $E_{l,k}\in\mathfrak{n}$ be the canonical basic element
$E_{l,k}(e_k)=e_l$ and $E_{l,k}(e_i)=0$ for $i\not=k$. The elements $E_{l,k}$, for $1\leq k<l\leq p+q+2$,
form a basis of $\mathfrak{n}$.

For $\underline{t}=(t_1,\ldots,t_7)\in\mathbb{C}^7$, we consider the element $x_{\underline{t}}\in\mathfrak{n}$ defined
as follows:

\smallskip
\begin{tabbing}
\qquad\qquad \= $x_{\underline{t}}(e_1)=t_1(e_2+t_2(e_{p+1}+t_3e_{p+2})),$ \\[1mm]
\> $x_{\underline{t}}(e_2)=t_2(e_{p+1}+t_3e_{p+2}), \quad x_{\underline{t}}(e_i)=0\ \mbox{for $i=3,\ldots,p-1$},$ \\[1mm]
\> $x_{\underline{t}}(e_p)=t_2t_4(e_{p+1}+t_3e_{p+2}),$ \\[1mm]
\> $x_{\underline{t}}(e_{p+1})=\left\{\begin{array}{l}
t_3(e_{p+2}+t_5(e_{p+2}+t_6(e_{p+q+1}-t_4e_{p+3})))\ \ \mbox{if $q>2$}, \\
t_3(e_{p+2}+t_5(e_{p+2}+t_6(e_{p+3}-t_4e_{p+4})))\ \ \mbox{if $q=2$}, \\
\end{array}\right.$ \\[1mm]
\> $x_{\underline{t}}(e_{p+2})=\left\{\begin{array}{l}
t_6(e_{p+q+1}-t_4e_{p+3})\ \ \mbox{if $q>2$}, \\
t_6(e_{p+3}-t_4e_{p+4})\ \ \mbox{if $q=2$}, \\
\end{array}\right.$ \\[1mm]
\> $x_{\underline{t}}(e_i)=0\ \mbox{for $i=p+3,\ldots,p+q$},$ \\[1mm]
\> $x_{\underline{t}}(e_{p+q+1})=t_7e_{p+q+2}\quad\mbox{and}\quad x_{\underline{t}}(e_{p+q+2})=0.$
\end{tabbing}

\smallskip
\noindent
The map $\mathbb{C}^7\rightarrow\mathfrak{n}$, $\underline{t}\mapsto x_{\underline{t}}$ is well defined
and algebraic.
Let $T_\pi$ be the standard tableau associated to $\pi$ in the sense of section \ref{duality},
then the component ${\mathcal K}_\pi^{\mathrm{BC}}$ is the closure of the subset ${\mathcal B}_u^{T_\pi}$
for the Zariski topology.
A straightforward computation shows that $x_{\underline{t}}\in\varphi^{-1}(\Omega\cap {\mathcal B}_u^{T_\pi})$
whenever $t_3t_5t_6t_7\not=0$ (if $q>2$) or $t_3t_5t_6(t_4+t_7)\not=0$ (if $q=2$).
Hence $\varphi(x_{\underline{t}})\in\Omega\cap {\mathcal K}^{\mathrm{BC}}_\pi$ for all $\underline{t}\in\mathbb{C}^7$.
In particular, $F_0=\varphi(0)\in{\mathcal K}_\pi^{\mathrm{BC}}$.

For $i\in\{1,\ldots,7\}$, let $\underline{t}^{(i)}=(t_1,\ldots,t_7)$ with $t_i=t$ and $t_j=0$ for $j\not=i$.
Then the curve $\{x_{\underline{t}^{(i)}}:t\in\mathbb{C}\}$
lies in $\varphi^{-1}(\Omega\cap {\mathcal K}^{\mathrm{BC}}_\pi)$ and its tangent vector at $t=0$
is an element of $T_{0}\,\varphi^{-1}(\Omega\cap {\mathcal K}^{\mathrm{BC}}_\pi)$, the tangent space of  $\varphi^{-1}(\Omega\cap {\mathcal K}^{\mathrm{BC}}_\pi)$ at $0$.
Making $i$ run over $\{1,\ldots,7\}$, we get:
\[E_{2,1}, E_{p+1,2}, E_{p+2,p+1}, E_{p+q+1,p+2}, E_{p+q+2,p+q+1}\in T_{0}\,\varphi^{-1}(\Omega\cap {\mathcal K}^{\mathrm{BC}}_\pi).\]
The curves $\{x_{\underline{t}}:\underline{t}=(0,t,0,1,0,0,0),\ t\in\mathbb{C}\}$,
$\{x_{\underline{t}}:\underline{t}=(0,0,0,1,0,t,0),\ t\in\mathbb{C}\}$
and $\{x_{\underline{t}}:\underline{t}=(0,0,1,0,-1,t,0),\ t\in\mathbb{C}\}$ also lie in
$\varphi^{-1}(\Omega\cap {\mathcal K}^{\mathrm{BC}}_\pi)$ and pass through $0$ for $t=0$.
Considering their tangent vectors at $t=0$, we derive
\[E_{p+1,p},E_{l,p+2},E_{p+q+1,p+1}\in T_{0}\,\varphi^{-1}(\Omega\cap {\mathcal K}^{\mathrm{BC}}_\pi)\]
(with $l=p+3$ for $q>2$ and $l=p+4$ for $q=2$).
Altogether, we get $8$ linearly independent tangent vectors.

If $q=2$, then we obtain
\[\dim T_{F_0}\,{\mathcal K}_\pi^{\mathrm{BC}}
=\dim T_{0}\,\varphi^{-1}(\Omega\cap {\mathcal K}^{\mathrm{BC}}_\pi)\geq 8>7=\dim{\mathcal K}_\pi^{\mathrm{BC}}\]
and we are done. If $q>2$, then we construct some additional elements in the tangent space.

Fix $k\in\{3,\ldots,q\}$. For $t\in\mathbb{C}$, there is a unique element $g_t\in Z_u$ such that
\begin{eqnarray*}
&& g_t(e_{p+q+1})=\left\{\begin{array}{l}
e_{p+q+1}+te_{p+k+1}\ \mbox{if $k<q$,}\\
e_{p+q+1}+te_{p+q+2}\ \mbox{if $k=q$,}
\end{array}
\right.\\
&& g_t(e_i)=e_i\ \mbox{for $i\in\{p+1,p+2,p+q+2\}$}.
\end{eqnarray*}
As $F_0\in{\mathcal K}_\pi^{\mathrm{BC}}$, we get $g_tF_0\in {\mathcal K}_\pi^{\mathrm{BC}}$ for all $t\in\mathbb{C}$.
In fact we see that $g_tF_0=\varphi(y_t)$ with $y_t\in\mathfrak{n}$ defined by
\[y_t=t(\delta_{k=p}E_{2,1}+(E_{p+k,p}+...+E_{p+3,p-k+3})+\delta_{k=q}E_{p+q+2,p+q+1}),\]
where $\delta_{i=j}=1$ for $i=j$ and $\delta_{i=j}=0$ otherwise. We infer that
\[E_{p+k,p}+...+E_{p+3,p-k+3}\in T_{0}\,\varphi^{-1}(\Omega\cap {\mathcal K}^{\mathrm{BC}}_\pi)\ \ \forall k\in\{3,\ldots,q\}.\]
These tangent vectors are linearly independent of the previous ones we constructed.
Altogether we have $q+6$ tangent vectors.
It results:
\[\dim T_{F_0}\,{\mathcal K}_\pi^{\mathrm{BC}}
=\dim T_{0}\,\varphi^{-1}(\Omega\cap {\mathcal K}^{\mathrm{BC}}_\pi)\geq q+6>q+5=\dim{\mathcal K}_\pi^{\mathrm{BC}}.\]
Therefore, $F_0$ is a singular point of the component ${\mathcal K}_\pi^{\mathrm{BC}}$.
\hfill $\Box$

\subsection{Singular components of the form ${\mathcal K}_\pi^{\mathrm{BC}}$ with $\pi=(2,p,2)$}

In \cite[\S 2.3]{Fresse-Melnikov-1}, another example of a singular Bala-Carter component is given,
and this is ${\mathcal K}_\pi^{\mathrm{BC}}$ for $\pi=(2,3,2)$.
More generally, we show:

\begin{lemma}
\label{lemma-singularity-2-p-2}
Let $p\geq 3$ and $\pi=(2,p,2)$. Then, the component ${\mathcal K}_\pi^{\mathrm{BC}}$ is singular.
\end{lemma}

\noindent
{\em Proof.}
The arguments are of the same type as in the proof of Lemma \ref{lemma-singularity-1-p-q-1}.
We may suppose that $p>3$, since otherwise $\pi=(2,3,2)$
and the singularity of ${\mathcal K}_\pi^{\mathrm{BC}}$ follows from \cite[\S 2.3]{Fresse-Melnikov-1}.
The considered component ${\mathcal K}_\pi^{\mathrm{BC}}$ lies in ${\mathcal B}_u$ for $u\in\mathrm{End}(\mathbb{C}^{p+4})$
nilpotent of Jordan form $\lambda(u)=(p,2,2)$.
Due to formula (\ref{relation-dimension-Bu}),
$\dim{\mathcal K}_\pi^{\mathrm{BC}}=\dim{\mathcal B}_u=6$.

We fix a Jordan basis $(e_1,\ldots,e_{p+4})$, on which $u$ acts as follows:
\[
\begin{array}{ccccccccccccccccccccccccccccccccc}
0 & \leftarrow & e_1 & \leftarrow & e_4 & \leftarrow & e_5 & \cdots & \leftarrow & e_{p+1} & \leftarrow & e_{p+4} \\
0 & \leftarrow & e_2 & \leftarrow & e_{p+2} \\
0 & \leftarrow & e_3 & \leftarrow & e_{p+3}
\end{array}
\]
We consider the adapted flag $F_0=(\langle e_1,\ldots,e_i\rangle)_{i=0,\ldots,p+4}$.
Our purpose is to show that $F_0$ is a singular point of the component ${\mathcal K}_\pi^{\mathrm{BC}}$.
To do this, as in the proof of Lemma \ref{lemma-singularity-1-p-q-1}, we show that $F_0$ lies in
${\mathcal K}_\pi^{\mathrm{BC}}$, and we construct more than $6$ linearly independent vectors
of the tangent space $T_{F_0}\,{\mathcal K}_\pi^{\mathrm{BC}}$.

As in the proof of Lemma \ref{lemma-singularity-1-p-q-1}, we consider the group
$B=\{g\in GL(V):ge_i\in\langle e_{i},\ldots,e_{p+4}\rangle\ \forall i\}$ and
the orbit $\Omega=B F_0$, which is an open affine neighborhood of $F_0$ in the variety of complete
flags of $\mathbb{C}^{p+4}$. We consider the subgroup $U\subset B$ of unipotent matrices,
and its Lie algebra $\mathfrak{n}=\{g-I:g\in U\}$. The map $\varphi:\mathfrak{n}\rightarrow\Omega$,
$x\mapsto (I+x)F_0$ is an isomorphism of algebraic varieties.
For $1\leq k<l\leq p+4$, let $E_{l,k}(e_k)=e_l$ and $E_{l,k}(e_i)=0$ for $i\not=k$,
then the elements $E_{l,k}$ (for all $1\leq k<l\leq p+4$) form the canonical basis
of $\mathfrak{n}$.

For $\underline{t}=(t_1,\ldots,t_6)\in\mathbb{C}^6$, we consider $x_{\underline{t}}\in\mathfrak{n}$
defined as follows:

\smallskip
\begin{tabbing}
\qquad\qquad \= $x_{\underline{t}}(e_1)=t_1(e_2+t_2e_3),$ \\[1mm]
\> $x_{\underline{t}}(e_2)=t_2e_3+t_3(e_3+t_4t_5(e_4+t_1(e_{p+2}+t_2e_{p+3}))),$ \\[1mm]
\> $x_{\underline{t}}(e_3)=t_4t_5(e_4+t_1(e_{p+2}+t_2e_{p+3})),$ \\[1mm]
\> $x_{\underline{t}}(e_i)=0\ \mbox{for $i=4,\ldots,p$},$ \\[1mm]
\> $x_{\underline{t}}(e_{p+1})=t_4(e_{p+2}+t_2e_{p+3}),$ \\[1mm]
\> $x_{\underline{t}}(e_{p+2})=t_2e_{p+3}+t_6(e_{p+3}-t_1t_5e_{p+4}),$ \\[1mm]
\> $x_{\underline{t}}(e_{p+3})=-t_1t_5e_{p+4},$ \\[1mm]
\> $x_{\underline{t}}(e_{p+4})=0.$
\end{tabbing}

\smallskip
\noindent
The map $\mathbb{C}^6\rightarrow \mathfrak{n}$, $\underline{t}\mapsto x_{\underline{t}}$ is well defined
and algebraic.
The component ${\mathcal K}_\pi^{\mathrm{BC}}$ is the closure of the subset ${\mathcal B}_u^{T_\pi}$,
where $T_\pi$ is the standard tableau associated to the sequence $\pi$
(see section \ref{duality}).
It is straightforward to see that $x_{\underline{t}}\in\varphi^{-1}(\Omega\cap{\mathcal B}_u^{T_\pi})$
whenever $t_1t_3t_4t_5t_6\not=0$. Hence $x_{\underline{t}}\in\varphi^{-1}(\Omega\cap{\mathcal K}_\pi^{\mathrm{BC}})$
for all $\underline{t}\in\mathbb{C}^6$.
In particular, $F_0=\varphi(0)\in{\mathcal K}_\pi^{\mathrm{BC}}$.

For $i\in\{1,\ldots,6\}$ and $t\in\mathbb{C}$,
let $\underline{t}^{(i)}=(t_1,\ldots,t_6)$ with $t_i=t$ and $t_j=0$ for $j\not=i$.
The curves $\{x_{\underline{t}^{(i)}}:t\in\mathbb{C}\}$
(for $i=1,\ldots,6$) lie in $\varphi^{-1}(\Omega\cap{\mathcal K}_\pi^{\mathrm{BC}})$,
and their tangent vectors at $t=0$ provide elements
\[E_{2,1}, E_{3,2}, E_{p+2,p+1}, E_{p+3,p+2}\in T_0\,\varphi^{-1}(\Omega\cap{\mathcal K}_\pi^{\mathrm{BC}}).\]
The curves $\{x_{\underline{t}}:\underline{t}=(t,-1,1,0,0,1),\ t\in\mathbb{C}\}$,
$\{x_{\underline{t}}=(0,-1,1,t,0,1),\ t\in\mathbb{C}\}$ and $\{x_{\underline{t}}:\underline{t}=(t,0,0,0,1,0),\ t\in\mathbb{C}\}$
also pass through $0$ at $t=0$. Considering their tangent vectors at $t=0$, we get:
\[E_{3,1},E_{p+3,p+1},E_{p+4,p+3}\in T_0\,\varphi^{-1}(\Omega\cap{\mathcal K}_\pi^{\mathrm{BC}}).\]
Altogether, we obtain that $7$ linearly independent elements lie in the tangent space
$T_0\,\varphi^{-1}(\Omega\cap{\mathcal K}_\pi^{\mathrm{BC}})$. It results:
\[\dim T_{F_0}\,{\mathcal K}_\pi^{\mathrm{BC}}=\dim T_0\,\varphi^{-1}(\Omega\cap{\mathcal K}_\pi^{\mathrm{BC}})
\geq 7>6=\dim{\mathcal K}_\pi^{\mathrm{BC}}.\]
Therefore, $F_0$ is a singular point of ${\mathcal K}_\pi^{\mathrm{BC}}$.
\hfill $\Box$

\subsection{Proof of Proposition \ref{proposition-theorem-1-if}}

Consider a sequence $\pi=(\pi_1,\ldots,\pi_r)\in\Lambda_u$ such that
\[\pi\geq (1,2,2,1)\quad\mbox{or}\quad \pi\geq (2,3,2).\]

First, assume $\pi\geq (1,2,2,1)$, that is, there are $1\leq i_1<i_2<i_3<i_4\leq r$
such that $\pi_{i_1},\pi_{i_4}\geq 1$, $\pi_{i_2},\pi_{i_3}\geq 2$.
Set $\hat\pi=(\pi_{i_1},\pi_{i_2},\pi_{i_3},\pi_{i_4})$ and
$\hat\pi'=(1,\pi_{i_2},\pi_{i_3},1)$.
Using Lemma \ref{lemma-singularity-1-p-q-1}, we get that the Bala-Carter component
${\mathcal K}_{\hat\pi'}^{\mathrm{BC}}$ associated to $\hat\pi'$ is singular.
Using Corollary \ref{corollary-extremal-points}, we infer that ${\mathcal K}_{\hat\pi}^{\mathrm{BC}}$ is singular.
Finally, applying Corollary \ref{corollary-connected-component},
it follows that the component ${\mathcal K}_{\pi}^{\mathrm{BC}}$ is singular.

Next, assume $\pi\geq (2,3,2)$. Hence, there are $1\leq i_1<i_2<i_3\leq r$
such that $\pi_{i_1},\pi_{i_3}\geq 2$, $\pi_{i_2}\geq 3$.
Similarly as above, we set $\hat\pi=(\pi_{i_1},\pi_{i_2},\pi_{i_3})$ and
$\hat\pi'=(2,\pi_{i_2},2)$. By Lemma \ref{lemma-singularity-2-p-2},
the component ${\mathcal K}_{\hat\pi'}^{\mathrm{BC}}$ associated to $\hat\pi'$ is singular.
Corollary \ref{corollary-extremal-points} then implies that ${\mathcal K}_{\hat\pi}^{\mathrm{BC}}$ is singular.
Finally, using Corollary \ref{corollary-connected-component},
we infer that the component ${\mathcal K}_{\pi}^{\mathrm{BC}}$ is singular.

In each case, we obtain that the component ${\mathcal K}_{\pi}^{\mathrm{BC}}$ is singular.
The proof of Proposition \ref{proposition-theorem-1-if} is then complete.

\section{Smooth Bala-Carter components}

\label{section-5}

As in the previous section, we consider a sequence
$\pi=(\pi_1,\ldots,\pi_r)\in\Lambda_u$, and
its corresponding Bala-Carter component ${\mathcal K}_\pi^{\mathrm{BC}}\subset{\mathcal B}_u$.
Our purpose is now to establish the implication $(\Rightarrow)$
of Theorem \ref{theorem-1}, namely:

\begin{proposition}
\label{proposition-theorem-1-only-if}
If ${\mathcal K}_\pi^{\mathrm{BC}}$ is singular,
then $\pi\geq(1,2,2,1)$ or $\pi\geq(2,3,2)$.
\end{proposition}

We prove the contraposal: suppose $\pi\not\geq(1,2,2,1)$ and $\pi\not\geq(2,3,2)$,
and let us show that ${\mathcal K}_\pi^{\mathrm{BC}}$ is smooth.
Of course, if $\pi$ is of the form $\pi=(\pi_1)$ (one term in the sequence),
we know that the component is smooth (it is a point).
More generally, if $\pi$ is of the form $\pi=(1,\ldots,1,p,1,\ldots,1)$ (the
component is of hook type) or $\pi=(p,q)$ (the component is of two-row type), then
we know that it is smooth (cf. \cite{Fung}, \cite{Vargas}).
Taking also into account Corollary \ref{corollary-symmetry},
the situations which remain to be analyzed are the following two ones:
\[\mbox{a)}\ \ \pi=(p,1,\ldots,1,2,1,\ldots,1,q),\mbox{ \ and \ b)\ \ }\pi=(1,\ldots,1,p,1,\ldots,1,q).\]
In both situations, $p,q\geq 2$ are arbitrary, the number of $1$'s is arbitrary (can be zero).
Consider situation a), and let us show that the component is smooth in this case.
First suppose $p=q=2$: applying \cite[Theorem 1.2]{Fresse-Melnikov-2},
we obtain that the component ${\mathcal K}_\pi^{\mathrm{BC}}$ is smooth in this case.
Next, applying Corollary \ref{corollary-extremal-points}, we derive
that ${\mathcal K}_\pi^{\mathrm{BC}}$ is smooth for $p,q\geq 2$ arbitrary.

It remains to show that ${\mathcal K}_\pi^{\mathrm{BC}}$ is smooth, whenever $\pi$ accords with situation b).
Again, applying Corollary \ref{corollary-extremal-points}, we see that it is sufficient
to treat the case where $p=q$. Applying Corollary \ref{corollary-connected-component},
we see that we may suppose in addition that the number of $1$'s on the left of $p$ coincides
with the number of $1$'s on the right and that this number is nonzero. Therefore, it remains to show the following

\begin{proposition}
\label{proposition-p-1-p-1}
Let $p\geq 2$ be arbitrary, let $m\geq 1$. Set
\[\pi=(\underbrace{1,\ldots,1}_{\mbox{\rm\scriptsize $m$ terms}},p,\underbrace{1,\ldots,1}_{\mbox{\rm\scriptsize $m$ terms}},p).\]
Then, the Bala-Carter component ${\mathcal K}_\pi^{\mathrm{BC}}$ is smooth.
\end{proposition}

The remainder of this section is devoted to the proof of Proposition \ref{proposition-p-1-p-1}.

\subsection{Setting}

\label{section-setting}

In \cite[\S 4.3--4.5]{Fresse-Melnikov-1}, it is proved that the component ${\mathcal K}_\pi^{\mathrm{BC}}$
is smooth, for $\pi=(1,p,p)$.
The proof of Proposition \ref{proposition-p-1-p-1} we give here will follow the same pattern.

In the remainder of the section, we consider $V=\mathbb{C}^{2p+2m}$ and $u\in\mathrm{End}(V)$ nilpotent,
of Jordan form $\lambda(u)=(p,p,1,\ldots,1)$ (with $p\geq 2$, $m\geq 1$).
Let $\pi=(1,\ldots,1,p,1,\ldots,1,p)$ be as in the statement of Proposition \ref{proposition-p-1-p-1},
and let ${\mathcal K}_\pi^{\mathrm{BC}}\subset{\mathcal B}_u$ be the corresponding Bala-Carter component.
In general, the irreducible components of ${\mathcal B}_u$ are parameterized by standard tableaux
and we denote by ${\mathcal K}^T$ the component associated to the tableau $T$ (cf. section \ref{Spaltenstein-construction}).
In particular, the component ${\mathcal K}_\pi^{\mathrm{BC}}$ coincides
with ${\mathcal K}^{T_\pi}$, where $T_\pi$ is the following standard tableau
\[T_\pi=\begin{array}{cccccccccccccccccc}
\mbox{\scriptsize $1$} & \mbox{\scriptsize $m{+}2$} & \mbox{\scriptsize $\cdots$} & \mbox{\scriptsize $m{+}p$} \\
\mbox{\scriptsize $\vdots$} & \mbox{\scriptsize $2m{+}p{+}2$} & \mbox{\scriptsize $\cdots$} & \mbox{\scriptsize $2m{+}2p$} \\
\mbox{\scriptsize $m{+}1$} \\
\mbox{\scriptsize $m{+}p{+}1$} \\
\mbox{\scriptsize $\vdots$} \\
\mbox{\scriptsize $2m{+}p{+}1$}
\end{array}\]
(cf. section \ref{duality}).
The component admits ${\mathcal B}_u^{T_\pi}$ as a dense open subset (cf. section \ref{Spaltenstein-construction}).
Due to formula (\ref{relation-dimension-Bu}),
\[\dim {\mathcal K}_\pi^{\mathrm{BC}}=\dim{\mathcal B}_u=\frac{(2m+2)(2m+1)}{2}+p-1=2m^2+3m+p.\]

\subsection{Special flags $F_{\underline{d}}$}

\label{flags-Fd}

In this subsection, the purpose is to produce a set of special elements
$F_{\underline{d}}\in{\mathcal B}_u$ with the property that it will be sufficient
to check the smoothness of the component ${\mathcal K}_\pi^{\mathrm{BC}}$ at these particular points.

We fix a basis $(e_1,\ldots,e_{2p+2m})$ of $V$, on which $u$ acts as follows:
\begin{eqnarray*}
&& e_1,\ldots,e_{2m+2}\in\ker u, \\
&& u(e_{2m+3})=e_1,\ u(e_{2m+4})=e_2\ \mbox{and}\ u(e_{2m+i})=e_{2m+i-2}\ \forall i\in\{5,\ldots,2p\}.
\end{eqnarray*}
We can also describe
the action of $u$ on the basis with the following tableau $\tau_0$:
\[\tau_0=\begin{array}{cccccccccccccccccc}
\mbox{\scriptsize $1$} & \mbox{\scriptsize $2m{+}3$} & \mbox{\scriptsize $\cdots$} & \mbox{\scriptsize $2m{+}2p{-}1$} \\
\mbox{\scriptsize $2$} & \mbox{\scriptsize $2m{+}4$} & \mbox{\scriptsize $\cdots$} & \mbox{\scriptsize $2m{+}2p$} \\
\mbox{\scriptsize $3$} \\
\mbox{\scriptsize $\vdots$} \\
\mbox{\scriptsize $2m{+}2$}
\end{array}.\]
We have $u(e_i)=0$ if $i$ lies in the first column of $\tau_0$, and $u(e_i)=e_j$, where $j$ is the number
on the left of $i$ in $\tau_0$, otherwise.

We consider tableaux which are obtained by permuting the entries of $\tau_0$.
For a permutation $\sigma\in\mathbf{S}_{2p+2m}$, we denote by $\sigma(\tau_0)$ the tableau obtained from $\tau_0$
after replacing each entry $i$ by $\sigma_i$. The tableau $\sigma(\tau_0)$ is said to be row-increasing
if its entries increase from left to right along the rows.
Let $\mathbf{S}_u=\{\sigma\in\mathbf{S}_{2p+2m}:\sigma^{-1}(\tau_0)\mbox{ is row-increasing}\}$.

For $\sigma\in\mathbf{S}_{2p+2m}$ we consider the flag
$F_\sigma=(\langle e_{\sigma_1},\ldots,e_{\sigma_i}\rangle)_{i=0,\ldots,n}$.
It is easy to see that $F_\sigma\in{\mathcal B}_u$ if and only if $\sigma\in\mathbf{S}_u$.
We have (cf. \cite[Lemma 1.3]{Fresse}):

\begin{lemma}
\label{lemma-special-flags}
An irreducible component ${\mathcal K}^T\subset{\mathcal B}_u$ is smooth if and only if
every flag of the form $F_\sigma$ (with $\sigma\in\mathbf{S}_u$) lying in ${\mathcal K}^T$ is a smooth point.
\end{lemma}

According to this lemma, to determine if a component is smooth,
it is sufficient to study the smoothness of the special points of the form $F_\sigma$.
In the case of the component ${\mathcal K}_\pi^{\mathrm{BC}}$ which we are concerned with,
the set of special points to consider can be reduced.

Let $D$ denote the set of uples $\underline{d}=(d_1<\ldots<d_{m+1})$
with $d_i\in\{m+2,\ldots,m+p+i\}$ for all $i$, and such that
$|\{p+m,\ldots,p+2m+1\}\cap\{d_1,\ldots,d_{m+1}\}|\geq m$.
For $\underline{d}\in D$, we define a flag $F_{\underline{d}}\in{\mathcal B}_u$ as follows.
Write $\{m+2,\ldots,2p+2m\}\setminus\{d_1,\ldots,d_{m+1}\}=\{i_3,\ldots,i_{2p}\}$ with $i_3<\ldots<i_{2p}$.
Let $\tau_{\underline{d}}$ be the following tableau:
\[\tau_{\underline{d}}=\begin{array}{cccccccccccccccccc}
\mbox{\scriptsize $1$} & \mbox{\scriptsize $i_3$} & \mbox{\scriptsize $\cdots$} & \mbox{\scriptsize $i_{2p-1}$} \\
\mbox{\scriptsize $\vdots$} & \mbox{\scriptsize $i_4$} & \mbox{\scriptsize $\cdots$} & \mbox{\scriptsize $i_{2p}$} \\
\mbox{\scriptsize $m{+}1$} \\
\mbox{\scriptsize $d_1$} \\
\mbox{\scriptsize $\vdots$} \\
\mbox{\scriptsize $d_{m+1}$}
\end{array}.\]
The tableau $\tau_{\underline{d}}$ is row-increasing, hence it can be written
$\tau_{\underline{d}}=\sigma_{\underline{d}}^{-1}(\tau_0)$ with $\sigma_{\underline{d}}\in\mathbf{S}_u$.
We define $F_{\underline{d}}=F_{\sigma_{\underline{d}}}$.
Then, we have:

\begin{lemma}
\label{lemma-Fd}
The component ${\mathcal K}_\pi^{\mathrm{BC}}$ is smooth if and only if every flag
$F_{\underline{d}}$ (with $\underline{d}\in D$)
such that $F_{\underline{d}}\in{\mathcal K}_\pi^{\mathrm{BC}}$ is a smooth point of ${\mathcal K}_\pi^{\mathrm{BC}}$.
\end{lemma}

\noindent
{\em Proof.} The implication $(\Rightarrow)$ is immediate.
To prove the second implication, we rely on Lemma \ref{lemma-special-flags}.
Recall that $Z_u=\{g\in GL(V):gug^{-1}=u\}$ is the stabilizer of $u$.
Note that, if $F,F'\in {\mathcal K}_\pi^{\mathrm{BC}}$ are such that $F'\in \overline{Z_uF}$ is a smooth point of ${\mathcal K}_\pi^{\mathrm{BC}}$, then $F$ is also a smooth point of ${\mathcal K}_\pi^{\mathrm{BC}}$.
Then, to complete our proof, it is sufficient to show that
the closure of the $Z_u$-orbit of every flag of the form $F_\sigma$ (with $\sigma\in\mathbf{S}_u$)
lying in ${\mathcal K}_\pi^{\mathrm{BC}}$ contains a flag of the form $F_{\underline{d}}$ (with $\underline{d}\in D$).

Let $g_t\in GL(V)$ ($t\in\mathbb{C}$) be defined by $g_t(e_2)=e_2+te_1$,
$g_t(e_{i})=e_{i}+te_{i-1}$ for $i\in\{2m+4,\ldots,2m+2p\}$ even,
and $g_t(e_l)=e_l$ otherwise.
We have $g_t\in Z_u$.
Notice that $\lim_{t\rightarrow\infty}g_tF_\sigma=F_{\sigma'}$
where $\sigma'\in\mathbf{S}_u$ is such that
\[ \sigma'^{-1}_i=\min\{\sigma^{-1}_i,\sigma^{-1}_{i+1}\},\ \
\sigma'^{-1}_{i+1}=\max\{\sigma^{-1}_i,\sigma^{-1}_{i+1}\}
\]
for all $i\in\{1,2m+3,\ldots,2m+2p-1\}$ odd, and
\[ \sigma'^{-1}_i=\sigma^{-1}_i\mbox{\ \ for $i\in\{3,\ldots,2m+2\}$.}\]

Let $g'_t\in GL(V)$ ($t\in\mathbb{C}$) be defined by $g'_t(e_{2m+3})=e_{2m+3}+te_2$,
$g'_t(e_{i})=e_{i}+te_{i-1}$ for $i\in\{2m+5,\ldots,2m+2p-1\}$ odd,
and $g'_t(e_l)=e_l$ otherwise.
Likewise, $g'_t\in Z_u$ and $\lim_{t\rightarrow\infty}g'_tF_\sigma=F_{\sigma'}$
for $\sigma'\in\mathbf{S}_u$ such that
\[\sigma'^{-1}_2=\min\{\sigma^{-1}_2,\sigma^{-1}_{2m+3}\},\ \
\sigma'^{-1}_{2m+3}=\max\{\sigma^{-1}_{2},\sigma^{-1}_{2m+3}\},\]
\[ \sigma'^{-1}_{i}=\min\{\sigma^{-1}_{i},\sigma^{-1}_{i+1}\},\ \
\sigma'^{-1}_{i+1}=\max\{\sigma^{-1}_{i},\sigma^{-1}_{i+1}\}\]
for all $i\in\{2m+4,\ldots,2m+2p-2\}$ even, and
\[\sigma'^{-1}_{i}=\sigma^{-1}_{i}\mbox{\ \ for $i\in\{1,3,\ldots,2m+2,2m+2p\}$.}\]

For $1\leq i<j\leq 2m+2$ such that $(i,j)\not=(1,2)$,
let $h^{(i,j)}_t=h_t\in GL(V)$ ($t\in\mathbb{C}$) be defined by $h_t(e_j)=e_j+te_i$ and $h_t(e_l)=e_l$ for $l\not=j$.
Then $h_t\in Z_u$, and $\lim_{t\rightarrow\infty}h_tF_\sigma=F_{\sigma'}$
with $\sigma'\in\mathbf{S}_u$ such that
\[\sigma'^{-1}_i=\min\{\sigma^{-1}_i,\sigma^{-1}_j\},\ \ \sigma'^{-1}_j=\max\{\sigma^{-1}_i,\sigma^{-1}_j\},
\ \ \mbox{and}\ \ \sigma'^{-1}_l=\sigma^{-1}_l\mbox{\ \ for $l\notin\{i,j\}$}.\]

We assume $F_\sigma\in{\mathcal K}_\pi^{\mathrm{BC}}$.
Combining the operations $(\lim_{t\rightarrow \infty}g_t)$, $(\lim_{t\rightarrow \infty}g'_t)$
and $(\lim_{t\rightarrow \infty}h^{(i,j)}_t)$, we transform $F_\sigma$
into a flag $F_{\sigma'}$ with $\sigma'$ such that
\[\sigma'^{-1}_2<\sigma'^{-1}_i<\sigma'^{-1}_{i+1}\mbox{\ \ for all $i\in\{2m+3,\ldots,2m+2p-1\}$},\]
\[\sigma'^{-1}_1<\sigma'^{-1}_2<\ldots<\sigma'^{-1}_{2m+2}.\]
In particular, $F_{\sigma'}\in\overline{Z_uF_\sigma}$. As we suppose $F_\sigma\in{\mathcal K}_\pi^{\mathrm{BC}}$,
we get $F_{\sigma'}\in {\mathcal K}_\pi^{\mathrm{BC}}$. Let us show that $\sigma'=\sigma_{\underline{d}}$ for some $\underline{d}\in D$,
which will complete the proof.

Denote $c_i=\sigma'^{-1}_i$ for $i\in\{1,\ldots,m+1\}$, $d_i=\sigma'^{-1}_i$ for $i\in\{m+2,\ldots,2m+2\}$,
and $i_j=\sigma'^{-1}_{2m+j}$ for $j\in\{3,\ldots,2p\}$.
The tableau $\sigma'^{-1}(\tau_0)$ is as follows:
\[\sigma'^{-1}(\tau_0)=\begin{array}{cccccccccccccccccc}
\mbox{\scriptsize $c_{1}$} & \mbox{\scriptsize $i_3$} & \mbox{\scriptsize $\cdots$} & \mbox{\scriptsize $i_{2p-1}$} \\
\mbox{\scriptsize $\vdots$} & \mbox{\scriptsize $i_4$} & \mbox{\scriptsize $\cdots$} & \mbox{\scriptsize $i_{2p}$} \\
\mbox{\scriptsize $c_{m{+}1}$} \\
\mbox{\scriptsize $d_1$} \\
\mbox{\scriptsize $\vdots$} \\
\mbox{\scriptsize $d_{m+1}$}
\end{array}.\]
Moreover, we have $c_1<\ldots<c_{m+1}<d_1<\ldots<d_{m+1}$ and $c_2<i_3<\ldots<i_{2p}$.

Following section \ref{Jordan-orbits}, the sequence $\pi=(\pi_1,\ldots,\pi_r)$
induces a standard partition of $\{1,\ldots,n\}$, which we also denote by $\pi$,
and the component ${\mathcal K}_\pi^{\mathrm{BC}}$ is the closure of the standard orbit ${\mathcal Z}_\pi$.
By definition of the orbit ${\mathcal Z}_\pi$,
any flag $F=(V_0,\ldots,V_{2m+2p})\in{\mathcal Z}_\pi$ satisfies
\begin{eqnarray*}
&& V_{m+1}\subset\ker u,\\
&& \dim V_{m+p+i}\cap\ker u\geq m+1+i\ \ \forall i\in\{1,\ldots,m+1\}, \\
&& \mbox{and}\ \ V_{2m+p+1}\subset u^{-1}(V_{m+p-1}).
\end{eqnarray*}
By lower semi-continuity of the map $F\mapsto \mathrm{rank}\,u_{|V_j/V_i}$,
these three properties are satisfied more generally
for all $F\in{\mathcal K}_\pi^{\mathrm{BC}}$.
In particular, $F_{\sigma'}=(V'_0,\ldots,V'_{2m+2p})$ satisfies them.
As $V'_{m+1}\subset\ker u$, using the definition of the flag $F_{\sigma'}$,
we get
\[|\{c_1,\ldots,c_{m+1}\}\cap\{1,\ldots,m+1\}|=\dim V'_{m+1}\cap \ker u=m+1,\]
hence $c_i=i$ for all $i$.
Similarly, for all $i\in\{1,\ldots,m+1\}$
we get
\[m+1+|\{d_1,\ldots,d_{m+1}\}\cap\{1,\ldots,m+p+i\}|=\dim V'_{m+p+i}\cap\ker u\geq m+1+i,\]
hence $d_i\leq m+p+i$.
Moreover, as $V'_{2m+p+1}\subset u^{-1}(V'_{m+p-1})$, any Jordan block of the induced map $u_{|V'_{2m+p+1}/V'_{m+p-1}}=0$
has size $1$, which implies that the entries $m+p,\ldots,2m+p+1$ lie in different rows of the tableau
$\sigma'^{-1}(\tau_0)$. In particular, at most $2$ among these $m+2$ numbers can lie in the first two rows
of the tableau, therefore $|\{d_1,\ldots,d_{m+1}\}\cap\{m+p,\ldots,2m+p+1\}|\geq m$.
It follows $\underline{d}\in D$, and finally $\sigma'=\sigma_{\underline{d}}$. The proof of the lemma is now complete.
\hfill $\Box$

\medskip

\subsection{Proof of Proposition \ref{proposition-p-1-p-1}}

\label{proof-proposition-p-1-p-1}

According to Lemma \ref{lemma-Fd}, to show that the component ${\mathcal K}_\pi^{\mathrm{BC}}$ is smooth,
it is sufficient to check the smoothness at the points of the form $F_{\underline{d}}$.
Thus, Proposition \ref{proposition-p-1-p-1} follows from the next proposition.

\begin{proposition}
\label{proposition-Fd}
For all $\underline{d}\in D$, we have $F_{\underline{d}}\in{\mathcal K}_\pi^{\mathrm{BC}}$,
and moreover $F_{\underline{d}}$ is a smooth point of ${\mathcal K}_\pi^{\mathrm{BC}}$.
\end{proposition}

Our purpose is then to establish Proposition \ref{proposition-Fd}.
To do this, we employ the same technique as in the proof of \cite[Proposition 4.3]{Fresse-Melnikov-1}.
Let us outline our proof.
For $\underline{d}\in D$, recall that $F_{\underline{d}}=(\langle e_{\sigma_{\underline{d}}(1)},\ldots,e_{\sigma_{\underline{d}}(i)}\rangle)_{i=0,\ldots,2m+2p}$
, where the basis $(e_1,\ldots,e_{2p+2m})$ and the permutation $\sigma_{\underline{d}}\in\mathbf{S}_{2p+2m}$ have
been introduced in section \ref{flags-Fd}.
We consider the Borel subgroup
\[B=\{g\in GL(V):g e_{\sigma_{\underline{d}}(i)}\in\langle e_{\sigma_{\underline{d}}(i)},\ldots,e_{\sigma_{\underline{d}}(2p+2m)}\rangle\ \ \forall i\}\]
and the orbit $\Omega_{\underline{d}}=B F_{\underline{d}}$.
The orbit $\Omega_{\underline{d}}$ is an open subset of the variety of complete flags of $V$,
and it is isomorphic to an affine space: for $F=(V_0,\ldots,V_{2p+2m})\in\Omega_{\underline{d}}$,
there is a unique basis $(\eta_1,\ldots,\eta_{2p+2m})$ such that $V_i=\langle \eta_1,\ldots,\eta_i\rangle$
and
\begin{equation}
\label{basis-omega}
\eta_i=e_{\sigma_{\underline{d}}(i)}+\sum_{j=i+1}^{2m+2p} \phi_{i,j}\,e_{\sigma_{\underline{d}}(j)}
\end{equation}
for some $\phi_{i,j}\in\mathbb{C}$. The maps $F\mapsto \phi_{i,j}$ are algebraic and the product map
$F\mapsto (\phi_{i,j})_{1\leq i<j\leq 2p+2m}$ is an isomorphism from $\Omega_{\underline{d}}$ to the affine space
$\mathbb{C}^{(p+m)(2p+2m-1)}$.
Then, we construct a closed immersion $\Phi_{\underline{d}}:\mathbb{C}^{2m^2+3m+p}\rightarrow \Omega_{\underline{d}}$
with the following two properties:

\smallskip
\noindent
\begin{tabbing}
\ \  \= (A)\ \ \= For $\underline{t}=(t_1,t_2,\ldots)\in\mathbb{C}^{2m^2+3m+p}$
with $t_i\not=0$ for all $i$, one has $\Phi_{\underline{d}}(\underline{t})\in{\mathcal B}_u^{T_\pi}$.\\[1mm]
\> (B) \> $\Phi_{\underline{d}}(0,\ldots,0)=F_{\underline{d}}$.
\end{tabbing}

\smallskip
\noindent
By (A), we get that $\Phi_{\underline{d}}$ is an isomorphism on a locally closed subset of the component
${\mathcal K}_\pi^{\mathrm{BC}}$. Since $\dim {\mathcal K}_\pi^{\mathrm{BC}}=2m^2+3m+p$ (see section \ref{section-setting}), the image of
$\Phi_{\underline{d}}$ is actually an open subset of ${\mathcal K}_\pi^{\mathrm{BC}}$.
By (B), we get that the flag $F_{\underline{d}}$ lies in the image of $\Phi_{\underline{d}}$,
hence it lies in ${\mathcal K}_\pi^{\mathrm{BC}}$ and is a smooth point.

To lead our construction, we need the preliminary construction summed up in the following lemma, which we quote from \cite[\S 4.4]{Fresse-Melnikov-1}.

\begin{lemma}
\label{lemma-vi}
For $\underline{t}=(t_3,\ldots,t_{p+1})\in\mathbb{C}^{p-1}$, there
are vectors $v_1,\ldots,v_{2p}$ depending algebraically on $\underline{t}$,
with the following properties: \\[1mm]
(a) $v_1=e_1$, $v_2=e_2$, $v_i-e_{2m+i}\in\langle e_{2m+i+1},\ldots,e_{2p+2m}\rangle$ for all $i\in\{3,\ldots,2p\}$. \\[1mm]
(b) $\langle v_1,\ldots,v_i\rangle$ is $u$-stable for all $i\in\{1,\ldots,2p\}$. \\[1mm]
(c) For $i\in\{3,\ldots,p+1\}$, $u(v_i)=v_{i-2}+t_i v_{i-1}$. In particular
$v_i\in\ker u^{i-1}$, and if $t_j\not=0$ for all $j$, then $v_i\notin\ker u^{i-2}$. \\[1mm]
(d) Noting $t'_1=t'_2=0$ and $t'_j=t'_{j-2}+t_j$ for $j\in\{3,\ldots,p+1\}$, then
\[v_i-e_{2m+i}-t'_ie_{2m+i+1}\in\langle e_{2m+i+2},\ldots,e_{2p+2m}\rangle\ \ \forall i\in\{3,\ldots,p+1\},\]
\[v_i-e_{2m+i}-t'_{2p-i+2}e_{2m+i+1}\in\langle e_{2m+i+2},\ldots,e_{2p+2m}\rangle\ \ \forall i\in\{p+2,\ldots,2p-1\}.\]
(e) If $\underline{t}=(0,\ldots,0)$, then $v_i=e_{2m+i}$ for all $i\in\{3,\ldots,2p\}$.
\end{lemma}

\noindent
{\em Proof.}
We consider the subspace generated by the basis vectors $e_1,e_2,$ $e_{2m+3},\ldots,$ $e_{2m+2p}$.
The action of $u$ on these vectors is represented by the following picture:
\[\begin{array}{lllllllllllll}
0 & \leftarrow & e_1 & \leftarrow & e_{2m+3} & \leftarrow & \cdots & \leftarrow & e_{2m+2p-1} \\
0 & \leftarrow & e_2 & \leftarrow & e_{2m+4} & \leftarrow & \cdots & \leftarrow & e_{2m+2p}\,.
\end{array}\]
In particular, note that $u$ restricts to an isomorphism from the subspace
$E_1:=\langle e_{2m+3},\ldots,e_{2m+2p}\rangle$ onto $E_2:=\langle e_1,e_2,e_{2m+3},\ldots,e_{2m+2p-2}\rangle$.
Let $\check{u}:E_2\rightarrow E_1$ be its inverse, that is $\check{u}(e_i)=e_{2m+i+2}$ for $i\in\{1,2\}$
and $\check{u}(e_i)=e_{i+2}$ for $i\in\{3,\ldots,2p-2\}$.
We then put $v_1=e_1$, $v_2=e_2$, and by induction $v_i=\check{u}(v_{i-2}+t_iv_{i-1})$ for all $i\in\{3,\ldots,p+1\}$.

For $i\in\{p+1,\ldots,2p\}$, we define vectors $(v_{j}^{(i)})_{1\leq j\leq i}$ in the following manner.
First, we set $v_j^{(p+1)}=v_j$, $\forall j\in\{1,\ldots,p+1\}$. Then for $i\geq p+2$, we set
$v_1^{(i)}=e_1$, $v_2^{(i)}=e_2$, and by induction: $v_j^{(i)}=\check{u}(v_{j-2}^{(i-1)})$ if $j\geq 3$.
For $i\in\{p+2,\ldots,2p\}$, we put $v_i=v_i^{(i)}$.

The fact that the vectors $v_i$, for $i=1,\ldots,2p$, are well defined
and satisfy properties (a)--(e) of the lemma then follows from \cite[Lemmas 4.3-4.4]{Fresse-Melnikov-1}.
\hfill $\Box$

\bigskip

We are ready to prove Proposition \ref{proposition-Fd}.

\medskip

\noindent
{\em Proof of Proposition \ref{proposition-Fd}.}

\smallskip

The constructions will rely on the Jordan basis $(e_1,\ldots,e_{2p+2m})$ introduced in section \ref{flags-Fd}.
We also consider the subspaces
$E=\langle e_1,e_2\rangle$, $F=\langle e_3,\ldots,e_{m+1}\rangle$ and
$G=\langle e_{m+2},\ldots,e_{2m+2}\rangle$.
Let $\mathfrak{n}_F\subset\mathrm{End}(F)$ be the subspace of strictly lower triangular maps,
i.e., $X\in\mathrm{End}(F)$ such that
$Xe_i\in\langle e_{i+1},\ldots,e_{m+1}\rangle$ for all $i=3,\ldots,m+1$.
Let $\mathfrak{n}_G\subset\mathrm{End}(G)$ be the subspace of strictly lower triangular maps.
Let ${\mathcal L}(E,F)$, ${\mathcal L}(F,G)$ be the spaces of linear maps
from $E$ to $F$ and $F$ to $G$ respectively.

Let $\underline{d}=(d_1,\ldots,d_{m+1})\in D$ and let $F_{\underline{d}}$ be the corresponding flag.
The purpose of the construction is to associate to a uple of variables $\underline{t}\in\mathbb{C}^{2m^2+3m+p}$
a certain basis $(\eta_1,\ldots,\eta_{2p+2m})$ of $V$ satisfying relation (\ref{basis-omega}),
such that the map
$\Phi_{\underline{d}}:\underline{t}\mapsto (\langle \eta_1,\ldots,\eta_i\rangle)_{i=0,\ldots,2m+2p}$
will fulfill the properties (A) and (B) above.
We distinguish two situations:

\smallskip
\noindent
\begin{tabbing}
\qquad \= 1)\ \ \= $d_1\geq p+m$, \\
\> 2)\ \ \> $d_1<p+m$.
\end{tabbing}

\smallskip
\noindent
We describe the uple of variables in each case.
In case (1), the uple of variables $\underline{t}$ is taken of the form
\[\underline{t}=(X,Y,H,K,s_1,t_3,\ldots,t_{p+1},y_1,y'_1,\ldots,y_{m+1},y'_{m+1}),\]
where $X\in \mathfrak{n}_F$, $Y\in \mathfrak{n}_G$,
$H\in{\mathcal L}(E,F)$, $K\in{\mathcal L}(F,G)$ and
$s_1,t_l,y_j,y'_j\in\mathbb{C}$.
In case (2), the uple $\underline{t}$ is taken of the form
\[\underline{t}=(X,Y,H,K,s_1,t_3,\ldots,t_{d_1{-}m{+}1},t_{d_1{-}m{+}3},\ldots,t_{p+2},y_1,y'_1,\ldots,y_{m+1},y'_{m+1}),\]
where $X,Y,H,K$ are as above, $s_1,t_l,y_j,y'_j\in\mathbb{C}$.
In fact, in all cases:
\[\underline{t}\in\mathfrak{n}_F\times\mathfrak{n}_G\times{\mathcal L}(E,F)\times{\mathcal L}(F,G)\times\mathbb{C}^{p+2(m+1)},\]
and it is straightforward to check that the space $\mathfrak{n}_F\times\mathfrak{n}_G\times{\mathcal L}(E,F)\times{\mathcal L}(F,G)\times\mathbb{C}^{p+2(m+1)}$ has dimension $2m^2+3m+p$, that is, the same dimension as the component ${\mathcal K}_\pi^{\mathrm{BC}}$.

In case (2), we set in addition $t_{d_1-m+2}=-y_1y'_1$. In both cases, we consider
the vectors $v_1,\ldots,v_{2p}$ associated to the uple $(t_3,\ldots,t_{p+1})$ in the sense of Lemma \ref{lemma-vi}.
We put $f_i=e_i+Xe_i+Ke_i$ for $i\in\{3,\ldots,m+1\}$ and $f_{d_j}=e_{m+1+j}+Ye_{m+1+j}$ for $j\in\{1,\ldots,m+1\}$.
We also put $f_1=e_1+s_1e_2+(I+X+K)He_1$ and $f_2=e_2+(I+X+K)He_2$, where $I$ is the identity endomorphism of $V$.
In both cases, put
\[\eta_i=f_i\ \ \mbox{for all $i\in\{3,\ldots,m+1\}$}.\]
We will define in each case the remaining vectors $\eta_1,\eta_2,\eta_{m+2},\ldots,\eta_{2p+2m}$.
It will follow from the definition that
$\eta_i-e_{\sigma_{\underline{d}}(i)}\in\langle e_{\sigma_{\underline{d}}(i+1)},\ldots,e_{\sigma_{\underline{d}}(2p+2m)}\rangle$ for all $i$,
and the so-obtained map $\Phi_{\underline{d}}:\mathbb{C}^{2m^2+3m+p}\rightarrow\Omega_{\underline{d}}$,
$\underline{t}\mapsto (\langle\eta_1,\ldots\eta_i\rangle)_{i=0,\ldots,2p+2m}$
will be algebraic.
In addition, we will verify the following four properties:

\smallskip
\noindent
\begin{tabbing}
\qquad \= a)\ \ \= $\Phi_{\underline{d}}(\mathbb{C}^{2m^2+3m+p})\subset{\mathcal B}_u$; \\[1mm]
\> b) \> the map $\Phi_{\underline{d}}$ is a closed immersion; \\[1mm]
\> c) \> $\Phi_{\underline{d}}(\underline{t})\in{\mathcal B}_u^{T_{\pi}}$
whenever all numbers $t_l,y_j,y'_j$ are nonzero; \\[1mm]
\> d) \> $\Phi_{\underline{d}}(0,\ldots,0)=F_{\underline{d}}$.
\end{tabbing}

\smallskip
\noindent
As a conclusion, $\Phi_{\underline{d}}(\mathbb{C}^{2m^2+3m+p})\cong \mathbb{C}^{2m^2+3m+p}$
is an open neighborhood of the flag $F_{\underline{d}}$ in the component ${\mathcal K}_\pi^{\mathrm{BC}}$,
therefore $F_{\underline{d}}$ is a smooth point of ${\mathcal K}_\pi^{\mathrm{BC}}$:
the proof will be complete.

\medskip
\noindent
(1) We suppose $d_1\geq p+m$.

\smallskip
By definition of the set $D$, the uple $\underline{d}$
satisfies $d_1<d_2<\ldots<d_{m+1}\leq 2m+p+1$.
Hence, there is a unique $k\in\{0,\ldots m+1\}$ such that
$p+m+k\notin\{d_1,\ldots,d_{m+1}\}$.
Consequently the indices $d_1,\ldots,d_{m+1}$ and the indices $i_3,\ldots,i_{2p}$
involved in the definition of $\sigma_{\underline{d}}$ are given by:
\[d_j=p+m+j-1\ \,\forall j\in\{1,\ldots,k\}\ \,\mbox{and}\ \,d_j=p+m+j\ \,\forall j\in\{k+1,\ldots,m+1\},\]
\[i_l=m-1+l\ \,\forall l\in\{3,\ldots,p\},\ \,i_{p+1}=p+m+k,\ \,\mbox{and \,}i_l=2m+l\ \,\forall l\in\{p+2,\ldots,2p\}.\]
The indices are organized as follows:
\[i_3<\ldots<i_p<d_1<\ldots<d_k<i_{p+1}<d_{k+1}<\ldots<d_{m+1}<i_{p+2}<\ldots<i_{2p}.\]
For $l\in\{1,\ldots,p\}$, $j\in\{1,\ldots,m+1\}$, we define a number $x_{l,j}$ by setting
$x_{p,j}=y_j$ and by induction $x_{l,j}=-t_{l+2}x_{l+1,j}$ for $l=2,\ldots,p-1$,
and $x_{1,j}=-(t_3-s_1)x_{2,j}$.
We put

\smallskip

\begin{tabbing}
\qquad \= $\eta_1=f_1+\textstyle{\sum\limits_{j=1}^{m+1}x_{1,j}f_{d_j}},\quad 
\eta_2=f_2+\textstyle{\sum\limits_{j=1}^{m+1}x_{2,j}f_{d_j}},$ \\[1mm]
\> $\eta_{i_l}=v_l+\textstyle{\sum\limits_{j=1}^{m+1}x_{l,j}f_{d_j}}\ \ \forall l\in\{3,\ldots,p\},$ \\[1mm]
\> $\eta_{d_j}=f_{d_j}+y'_jv_{p+1}\ \ \forall j\in\{1,\ldots,k\},$ \\[1mm]
\> $\eta_{i_{p+1}}=v_{p+1}+\textstyle{\sum\limits_{j=k+1}^{m+1}y'_jf_{d_j}},$ \\[1mm]
\> $\eta_{d_j}=f_{d_j}\ \ \forall j\in\{k+1,\ldots,m+1\},
\ \ \mbox{and}\ \ \eta_{i_l}=v_l\ \ \forall l\in\{p+2,\ldots,2p\}.$
\end{tabbing}

\smallskip
\noindent
Let us show properties a) -- d).

\smallskip
\noindent
a) Note that for $i\in\{1,\ldots,m+1,d_{k+1},\ldots,d_{m+1}\}$,
we have $\eta_i\in\ker u$, which implies $u(\eta_i)\in\langle \eta_1,\ldots,\eta_i\rangle$.
Using Lemma \ref{lemma-vi}\,(c) and the definition of $x_{l,j}$,
we have:
\[u(v_3)=v_1+t_3v_2=(e_1+s_1e_2+\textstyle{\sum\limits_{j=1}^{m+1}x_{1,j}f_{d_j}})
+(t_3-s_1)(e_2+\textstyle{\sum\limits_{j=1}^{m+1}x_{2,j}f_{d_j}}),\]
hence $u(\eta_{i_3})=u(v_3)\in\langle \eta_1,\ldots,\eta_{m+1}\rangle$. Similarly, for $l\in\{4,\ldots,p+1\}$:
\[u(v_l)=v_{l-2}+t_lv_{l-1}=(v_{l-2}+\textstyle{\sum\limits_{j=1}^{m+1}x_{l-2,j}f_{d_j}})
+t_l(v_{l-1}+\textstyle{\sum\limits_{j=1}^{m+1}x_{l-1,j}f_{d_j}}),\]
hence $u(\eta_{i_l})=u(v_l)\in\langle \eta_1,\ldots,\eta_{i_l}\rangle$ for
$l\in\{4,\ldots,p\}$, and $u(\eta_i)\in\mathbb{C}u(v_{p+1})\subset \langle \eta_1,\ldots,\eta_i\rangle$
for $i\in\{d_1,\ldots,d_k,i_{p+1}\}$.
Notice that $v_1,\ldots,v_{p+1}\in\langle\eta_1,\ldots,\eta_{2m+p+1}\rangle$.
Then, for $l\in\{p+2,\ldots,2p\}$, we get $u(\eta_{i_l})=v_{l-2}+t_lv_{l-1}\in\langle\eta_1,\ldots,\eta_{i_l}\rangle$.
Therefore, we have $u(\eta_i)\in\langle\eta_1,\ldots,\eta_i\rangle$ for all $i$. It follows
$\Phi_{\underline{d}}(\underline{t})\in{\mathcal B}_u$.

\smallskip
\noindent
b) We show that $\Phi_{\underline{d}}$ is a closed immersion, and to do this,
we show that the algebra morphism
$\Phi_{\underline{d}}^*:\mathbb{C}[\phi_{i,j}:1\leq i<j\leq 2p+2m]\rightarrow \mathbb{C}[X,Y,H,K,s_1,t_i,y_j,y'_j]$
associated to $\Phi_{\underline{d}}$ is surjective. The functions $\phi_{i,j}$ are those involved
in the expression of the basis $\eta_1,\ldots,\eta_{2p+2m}$ according to formula (\ref{basis-omega}).
Then, we see that the coefficients of the matrices $H,X,K,Y$ are recovered by considering
$\phi_{i,j}$ for $i\in\{1,2\}$, $j\in\{3,\ldots,m+1\}$, or
$i,j\in\{3,\ldots,m+1\}$, or $i\in\{3,\ldots,m+1\}$, $j\in\{d_1,\ldots,d_{m+1}\}$,
or $i,j\in\{d_1,\ldots,d_{m+1}\}$, respectively. Therefore, they belong to $\mathrm{Im}\,\Phi_{\underline{d}}^*$.
We have $s_1=\phi_{1,2}\in \mathrm{Im}\,\Phi_{\underline{d}}^*$.
By Lemma \ref{lemma-vi}\,(d), we have $t'_l=\phi_{i_l,i_{l+1}}\in \mathrm{Im}\,\Phi_{\underline{d}}^*$
for all $l\in\{3,\ldots,p+1\}$. Due to the definition of $t'_l$ given in Lemma \ref{lemma-vi},
we infer that $t_l\in\mathrm{Im}\,\Phi_{\underline{d}}^*$ for all $l\in\{3,\ldots,p+1\}$.
Considering $\phi_{i_{p},d_j}$, we see that $y_j\in \mathrm{Im}\,\Phi_{\underline{d}}^*$ for all $j$.
We have $y'_j=\phi_{d_j,i_{p+1}}\in\mathrm{Im}\,\Phi_{\underline{d}}^*$ for $j\leq k$, and
moreover considering $\phi_{i_{p+1},d_j}$
we see that $y'_j\in\mathrm{Im}\,\Phi_{\underline{d}}^*$
for all $j>k$.
Finally, we obtain that $\Phi_{\underline{d}}^*$ is surjective.

\smallskip
\noindent
c) Suppose that all numbers $t_3,\ldots,t_{p+1},y'_1$ are nonzero, and let us show that
$\Phi_{\underline{d}}(\underline{t})\in{\mathcal B}_u^{T_\pi}$.
Let $u_i$ denote the restriction of $u$ to the subspace $\langle \eta_1,\ldots,\eta_i\rangle$,
and let $\lambda(u_i)$ be the sequence of its Jordan block sizes.
First, notice that $\eta_1,\ldots,\eta_{m+1}\in\ker u$, hence for all $i\in\{1,\ldots,m+1\}$,
$\lambda(u_i)=(1,\ldots,1)$.
Next, we have either $u^{p-1}(\eta_{m+p})=y'_1u^{p-1}(v_{p+1})$ (if $k>0$)
or $u^{p-1}(\eta_{m+p})=u^{p-1}(v_{p+1})$ (if $k=0$). By Lemma \ref{lemma-vi}\,(c),
$u^{p-1}(v_{p+1})\not=0$.
It follows $\lambda(u_{m+p})=(p,1,\ldots,1)$ and necessarily
$\lambda(u_i)=(i-m,1,\ldots,1)$ for all $i\in\{m+2,\ldots,m+p\}$.
Also, as $f_{d_1},\ldots,f_{d_{m+1}}\in\ker u$, we have
$\lambda(u_i)=(p,1,\ldots,\ldots,1)$ for all $i\in\{m+p+1,\ldots,2m+p+1\}$.
Finally, as $\lambda(u_{2m+2p})=\lambda(u)=(p,p,1,\ldots,1)$, we derive
$\lambda(u_i)=(p,i-2m-p,1,\ldots,1)$ for all $i\in\{2m+p+2,\ldots,2m+2p\}$.
For all $i$, the sequence $\lambda(u_i)$ indeed coincides with the lengths of the rows of the subtableau of $T_\pi$
of entries $1,\ldots,i$. Therefore, $\Phi_{\underline{d}}(\underline{t})\in{\mathcal B}_u^{T_\pi}$.

\smallskip
\noindent
d) follows from Lemma \ref{lemma-vi}\,(e)
and the definition of the vectors $\eta_i$.

\medskip
\noindent
(2) We suppose $d_1< p+m$.

\smallskip

We always have $d_1\geq m+2$. In particular, in this case, we have $p\geq 3$.
Due to the definition of $D$,
we have $\{d_2,\ldots,d_{m+1}\}\subset\{p+m,\ldots,p+2m+1\}$,
so that there are exactly two numbers $h,k\in\{0,\ldots,m+1\}$ with $h<k$ such that
$p+m+h,p+m+k\notin\{d_1,\ldots,d_{m+1}\}$. Then, the indices $d_j$ and $i_l$ are given by
\[d_j=p+m+j-2\ \ \forall j\in\{2,\ldots,h+1\},\quad
d_j=p+m+j-1\ \ \forall j\in\{h+2,\ldots,k\},\]
\[d_j=p+m+j\ \ \forall j\in\{k+1,\ldots,m+1\},\]
and, letting $d_1=m+c$,
\[i_l=m-1+l\ \ \forall l\in\{3,\ldots,c\},\quad
i_l=m+l\ \ \forall l\in\{c+1,\ldots,p-1\},\]
\[i_p=m+p+h,\quad i_{p+1}=m+p+k,\quad i_l=2m+l\ \ \forall l\in\{p+2,\ldots,2p\}.\]
The indices are then organized as follows:
\[i_3<\ldots<i_c<d_1<i_{c+1}<\ldots<i_{p-1}<d_2<\ldots<d_{h+1}<i_p<d_{h+2}<\ldots\]
\[<d_k<i_{p+1}<d_{k+1}<\ldots<d_{m+1}<i_{p+2}<\ldots<i_{2p}.\]
For $l\in\{1,\ldots,c\}$, we define a number $x_l$ by setting $x_c=y_1$,
and by induction $x_l=-t_{l+2}x_{l+1}$ for $l\in\{2,\ldots,c-1\}$, and $x_1=-(t_3-s_1)x_2$.
Then, we put:

\smallskip

\begin{tabbing}
\qquad \= $\eta_1=f_1+x_1f_{d_1},\quad \eta_2=f_2+x_2f_{d_1},$ \\[1mm]
\> $\eta_{i_l}=v_l+x_l f_{d_1}\ \ \forall l\in\{3,\ldots,c\},\quad \eta_{d_1}=f_{d_1}+y'_1v_{c+1},$ \\[1mm]
\> $\eta_{i_l}=v_l+t_{l+2}v_{l+1}\ \ \forall l\in\{c+1,\ldots,p-1\},$ \\[1mm]
\> $\eta_{d_j}=f_{d_j}+y_jv_p+y'_jv_{p+1}\ \ \forall j\in\{2,\ldots,h+1\},$ \\[1mm]
\> $\eta_{i_p}=v_p+t_{p+2}v_{p+1}+\textstyle{\sum\limits_{j=h+2}^{m+1}y_jf_{d_j}},$ \\[1mm]
\> $\eta_{d_j}=f_{d_j}+y'_jv_{p+1}\ \ \forall j\in\{h+2,\ldots,k\},\quad
\eta_{i_{p+1}}=v_{p+1}+\textstyle{\sum\limits_{j=k+1}^{m+1}y'_jf_{d_j}},$ \\[1mm]
\> $\eta_{d_j}=f_{d_j}\ \ \forall j\in\{k+1,\ldots,m+1\},\ \ \mbox{and}\ \ \eta_{i_l}=v_l\ \ \forall l\in\{p+2,\ldots,2p\}.$
\end{tabbing}

\smallskip
\noindent
We show properties a) -- d).

\smallskip
\noindent
a) Let us check that $u(\eta_i)\in\langle \eta_1,\ldots,\eta_i\rangle$ for all $i$.
First, notice that $\eta_1,\ldots,\eta_{m+1}\in\ker u$.
Similarly as in case (1), using Lemma \ref{lemma-vi}\,(c), we have:
\[u(v_3)=v_1+t_3v_2=(e_1+s_1e_2+x_1f_{d_1})+(t_3-s_1)(e_2+x_2f_{d_1})\in\langle \eta_1,\ldots,\eta_{m+1}\rangle,\]
and for $l\in\{4,\ldots,c+1\}$:
\[u(v_l)=v_{l-2}+t_lv_{l-1}=(v_{l-2}+x_{l-2}f_{d_1})+t_l(v_{l-1}+x_{l-1}f_{d_1})\in\langle \eta_i:1\leq i\leq i_{l-1}\rangle.\]
Using that $x_c=y_1$ and $t_{c+2}=t_{d_1-m+2}=-y_1y'_1$, we get:
\[u(v_{c+2})=v_c+t_{c+2}v_{c+1}=(v_c+x_cf_{d_1})-y_1(f_{d_1}+y'_1v_{c+1})\in\langle\eta_i:1\leq i\leq d_1\rangle.\]
Moreover, $u(v_l)=v_{l-2}+t_lv_{l-1}=\eta_{i_{l-2}}$ for $l\in\{c+3,\ldots,p+1\}$.
Then,
using in addition that $u(f_{d_j})=0$ for all $j$,
we get: $u(\eta_i)\in\langle\eta_1,\ldots,\eta_i\rangle$ for all $i\in\{d_1,\ldots,d_{m+1},i_3,\ldots,i_{p+1}\}$.
Observe that $v_1,\ldots,v_{p+1}\in\langle \eta_1,\ldots,\eta_{i_{p+2}-1}\rangle$.
Applying Lemma \ref{lemma-vi}\,(b),
we derive that $u(\eta_i)\in\langle\eta_1,\ldots,\eta_i\rangle$ also for $i\in\{i_{p+2},\ldots,i_{2p}\}$,
which completes the proof of property a).

\smallskip
\noindent
b) To show that the map $\Phi_{\underline{d}}$ is a closed immersion,
as in case (1), we show that its associated algebra morphism $\Phi_{\underline{d}}^*$ is surjective.
Exactly as in case (1), we obtain that the coefficients of the matrices $X,Y,H,K$ lie in the image of
$\Phi_{\underline{d}}^*$. We have $s_1=\phi_{1,2}\in\mathrm{Im}\,\Phi_{\underline{d}}^*$.
Let $t'_l$ be the numbers introduced in Lemma \ref{lemma-vi}\,(d).
For $l\in\{3,\ldots,p\}$, we get
$t'_l=t'_{2p+2-(2p-l+2)}=\phi_{i_{2p-l+2},i_{2p-l+3}}$
by applying Lemma \ref{lemma-vi}\,(d).
Similarly, we have $t'_{p+1}=\phi_{i_{p+1},i_{p+2}}$
and $t'_p+t_{p+2}=\phi_{i_p,i_{p+1}}$.
Then, in view of the definition of the numbers $t'_l$,
we infer that $t_l\in\mathrm{Im}\,\Phi_{\underline{d}}^*$
for all $l\in\{3,\ldots,p+2\}$.
We see that $x_l\in\mathrm{Im}\,\Phi_{\underline{d}}^*$ for all $l\in\{1,\ldots,c\}$,
by considering $\phi_{1,d_1}$, $\phi_{2,d_1}$ and $\phi_{i_l,d_1}$ for $l\geq 3$. 
Thus in particular $y_1=x_c\in \mathrm{Im}\,\Phi_{\underline{d}}^*$.
Also, $y'_1=\phi_{d_1,i_{c+1}}\in \mathrm{Im}\,\Phi_{\underline{d}}^*$.
We have $y_j=\phi_{d_j,i_p}\in\mathrm{Im}\,\Phi_{\underline{d}}^*$ for all $j\in\{2,\ldots,h+1\}$,
while for $j> h+1$ we see that $y_j\in\mathrm{Im}\,\Phi_{\underline{d}}^*$ by considering $\phi_{i_p,d_j}$.
Similarly, for $j\in\{2,\ldots,m+1\}$, we get that $y'_j\in\mathrm{Im}\,\Phi_{\underline{d}}^*$
by considering either $\phi_{d_j,i_{p+1}}$ (if $j\leq k$) or $\phi_{i_{p+1},d_j}$ (for $j>k$).
Finally, we have shown that $\Phi_{\underline{d}}^*$ is surjective.

\smallskip
\noindent
c) Suppose that all numbers $t_3,\ldots,t_{p+2},y_1,y'_1,y'_2$ are nonzero, and let us show that
$\Phi(\underline{t})\in{\mathcal B}_u^{T_\pi}$. As in case (1), we denote by $u_i$ the restriction
of $u$ to the subspace $\langle\eta_1,\ldots,\eta_i\rangle$, and by $\lambda(u_i)$ its Jordan form.
As in case (1), since $\eta_1,\ldots,\eta_{m+1}\in\ker u$, we have $\lambda(u_i)=(1,\ldots,1)$ for all $i\in\{1,\ldots,m+1\}$.
We have $u^{p-1}(v_{p+1})\not=0$ by Lemma \ref{lemma-vi}\,(c), which implies $u^{p-1}(\eta_{m+p})\not=0$,
therefore $\lambda(u_i)=(i-m,1,\ldots,1)$ for all $i\in\{m+2,\ldots,m+p\}$.
Notice that $\ker u\subset \langle\eta_1,\ldots,\eta_{2m+p+1}\rangle$,
and since $\dim\ker u=2m+2$, we have necessarily
$\lambda(u_{2m+p+1})=(p,1,\ldots,\ldots,1)$ (with $2m+1$ terms $1$) and consequently
$\lambda(u_i)$ is of the form $(p,1,\ldots,1)$ for all $i\in\{m+2,\ldots,2m+p+1\}$.
Finally, as in case (1), since $\lambda(u_{2p+2m})=\lambda(u)=(p,p,1,\ldots,1)$,
we derive that $\lambda(u_i)=(p,i-2m-p,1,\ldots,1)$
for all $i\in\{2m+p+2,\ldots,2p+2m\}$. For all $i$, the sequence $\lambda(u_i)$ coincides
with the lengths of the rows of the subtableau of $T_\pi$ of entries $1,\ldots,i$.
Therefore, $\Phi_{\underline{d}}(\underline{t})\in{\mathcal B}_u^{T_\pi}$.

\smallskip
\noindent
d) follows from Lemma \ref{lemma-vi}\,(e) and the definition of the vectors $\eta_1,\ldots,\eta_{2p+2m}$.

\medskip
The proof of Proposition \ref{proposition-Fd} is then complete.
\hfill $\Box$

\section{Components as iterated fiber bundles over projective spaces}

\label{section-6}

Recall that, for $T$ a standard tableau, we denote by $T^*$ its transpose
(that is, the rows of $T^*$ coincide with the columns of $T$).
To the tableaux $T$ and $T^*$, we associate the components ${\mathcal K}^T$ and ${\mathcal K}^{T^*}$,
imbedded in the appropriate Springer fibers ${\mathcal B}_u$ and ${\mathcal B}_{u^*}$.
In this section, we consider the situation where ${\mathcal K}^{T}$ contains a dense Jordan orbit,
and try to derive properties for ${\mathcal K}^{T^*}$.
The purpose of the section is to show Theorem \ref{theorem-3}, which states that in this situation,
${\mathcal K}^{T^*}$ is an iterated fiber bundle of base a sequence of projective spaces.
The proof is done by induction, and uses the combinatorial description
of components with a dense Jordan orbit, provided by Proposition \ref{proposition-2}.

\subsection{Concatenation of standard tableaux}

\label{sum-tableaux}

In this subsection, we show a preliminary result,
saying that if the standard tableau $T$ is obtained as the concatenation (which we will call the sum) of two standard tableaux $T_1,T_2$,
then the component ${\mathcal K}^T$ is isomorphic to the product of the components ${\mathcal K}^{T_1}$, ${\mathcal K}^{T_2}$.
We start with the definition of the sum of two standard tableaux.

Let $Y_1,Y_2$ be two Young diagrams with $n_1$ and $n_2$ boxes, respectively.
We define the sum $Y=Y_1+Y_2$ as the Young diagram with $n_1+n_2$ boxes
such that for all $j$, the length of the $j$-th row of $Y$ is the sum of the lengths
of the $j$-th rows of $Y_1$ and $Y_2$.
Now, let $T_1,T_2$ be two standard tableaux of shape $Y_1$ and $Y_2$ respectively.
We define the sum $T=T_1+T_2$ as the standard tableau of shape $Y$ such that for all $j$,
the $j$-th row of $T$ contains the entries of the $j$-th row of $T_1$, and the entries of the $j$-th row of $T_2$
increased by $n_1$. For instance:
\[\mbox{If}\ \ T_1=\young(12,3,4,5)\mbox{ \ \ and \ \ }T_2=\young(13,24)\,,\ \ \ \mbox{then\ \ }T=T_1+T_2=\young(1268,379,4,5)\,.\]

To the tableaux $T_1,T_2,T$, we associate respective components ${\mathcal K}^{T_1},{\mathcal K}^{T_2},{\mathcal K}^T$
imbedded in the appropriate Springer fibers.

\begin{proposition}
\label{proposition-somme-tableaux}
Suppose $T=T_1+T_2$. Then, ${\mathcal K}^T\cong {\mathcal K}^{T_1}\times{\mathcal K}^{T_2}$.
\end{proposition}

\noindent
{\em Proof.} Let $V=\mathbb{C}^{n_1+n_2}$, and let $u\in\mathrm{End}(V)$ be nilpotent
of Jordan form $Y(u)=Y$, so that the component ${\mathcal K}^T$ is imbedded in the Springer fiber ${\mathcal B}_u$.
For a $u$-stable subspace $W\subset V$, we denote by $Y(u_{|W})$ and $Y(u_{|V/W})$ the Young diagrams representing
the Jordan forms of the maps induced by $u$ on $W$ and $V/W$ respectively.
Let ${\mathcal G}_m(V)$ be the set of $m$-dimensional $u$-stable subspaces of $V$.
Our first step is to show that the set
\[{\mathcal A}=\{W\in{\mathcal G}_{n_1}(V):Y(u_{|W})=Y_1\mbox{ \ and \ }Y(u_{|V/W})=Y_2\}\]
is a single point.

\smallskip

We show this by induction on $n_1+n_2$ with immediate initialization if $n_1+n_2=0$.
Assume $n_1+n_2>0$. As in the proof of Proposition \ref{proposition-symmetry}, we denote by $\tilde V$ the space of linear forms $\phi:V\rightarrow \mathbb{C}$
and by $\tilde u:\tilde V\rightarrow \tilde V$ the dual nilpotent map of $u$.
If $W\subset V$ is a subspace, let $W^\perp=\{\phi\in\tilde V:W\subset\ker \phi\}$ be its dual space.
Let
$\tilde{\mathcal A}=\{W\in{\mathcal G}_{n_2}(\tilde V):Y(\tilde u_{|W})=Y_2\mbox{ \ and \ }Y(\tilde u_{|\tilde V/W})=Y_1\}$.
Notice that $W\in{\mathcal A}$ if and only if $W^\perp\in\tilde{\mathcal A}$.
Therefore, up to considering $\tilde{\mathcal A}$ instead of ${\mathcal A}$, we may assume
that the length of the first column of $Y_1$ is bigger than or equal to the length of the first column of $Y_2$.
Thus, the lengths of the first columns of $Y_1$ and $Y$ are equal,
which implies that every $W\in{\mathcal A}$ satisfies $\dim \ker (u_{|W})=\dim \ker u$,
hence $\ker u\subset W$. Let $r=\dim \ker u$.
Let $Y'_1$ be the Young diagram obtained from $Y_1$ by deleting the first column.
Let $V'=V/\ker u$ and let $u'\in\mathrm{End}(V')$ be the nilpotent map induced by $u$.
Thus $Y(u')=Y'_1+Y_2$.
Let ${\mathcal A}'=\{W\in{\mathcal G}_{n_1-r}(V'):Y(u'_{|W})=Y'_1\mbox{ \ and \ }Y(u'_{|V'/W})=Y_2\}$.
It is easy to see that $W\in{\mathcal A}$ if and only if $W/\ker u\in{\mathcal A}'$.
By induction hypothesis, ${\mathcal A}'$ is single point. It results that ${\mathcal A}$ is a single point,
as it was claimed.

\smallskip

Denote by $W_1$ the unique element of ${\mathcal A}$.

\smallskip

For $l\in\{1,2\}$ and $j\in\{1,\ldots,n_l\}$, let $Y_j^{T_l}$ be the shape of the subtableau of $T_l$ of entries
$1,\ldots,j$. Alternatively, for $0\leq i<j\leq n_1+n_2$, let $Y_{j/i}^T$ be the shape of
the rectification by jeu de taquin of the skew subtableau of $T$ of entries $i+1,\ldots,j$
(we refer to \cite{Fulton} for the definition of jeu de taquin).
By \cite[Theorem 3.3]{vanLeeuwen}, the set $\{F=(V_0,\ldots,V_{n_1+n_2})\in{\mathcal K}^T:Y(u_{|V_j/V_i})=Y_{j/i}^T\}$
is a nonempty open subset of the component ${\mathcal K}^T$.
On the other hand, by definition of $T$,
we have
$Y^T_{j/0}=Y_{j}^{T_1}$ for all $j\in\{1,\ldots,n_1\}$ and
$Y^T_{j/n_1}=Y_{j-n_1}^{T_2}$ for all $j\in\{n_1+1,\ldots,n_1+n_2\}$.
It follows that $\{F\in{\mathcal K}^T:V_{n_1}\in{\mathcal A}\}$ is a nonempty open subset of ${\mathcal K}^T$,
and we infer that $V_{n_1}=W_1$ for all $(V_0,\ldots,V_{n_1+n_2})\in{\mathcal K}^T$.

\smallskip

Let $u_1=u_{|W_1}$ and $u_2=u_{|V/W_1}$ be the nilpotent maps induced by $u$.
Let ${\mathcal B}_1$ (resp. ${\mathcal B}_2$) be the variety of complete flags of $W_1$ (resp. of $V/W_1$),
and let ${\mathcal B}_{u_1}\subset{\mathcal B}_1$ (resp. ${\mathcal B}_{u_2}\subset{\mathcal B}_2$)
be the subvariety of $u_1$-stable (resp. $u_2$-stable) flags. Thus the components ${\mathcal K}^{T_1},{\mathcal K}^{T_2}$
are imbedded in ${\mathcal B}_{u_1}$ and ${\mathcal B}_{u_2}$ respectively.
The map
\[\Phi:\{F\in{\mathcal B}_u:V_{n_1}=W_1\}\rightarrow {\mathcal B}_{u_1}\times{\mathcal B}_{u_2},\ \
(V_0,\ldots,V_{n_1+n_2})\mapsto ((V_j)_{j=0}^{n_1},(V_j/W_1)_{j=n_1}^{n_2})\]
is an isomorphism of algebraic varieties. Hence, it induces an isomorphism of ${\mathcal K}^T$ onto its image.
Let $\Phi_1,\Phi_2$ be the projections of $\Phi$ on both terms. We have
\[\{F\in{\mathcal K}^T:Y(u_{|V_j})=Y_j^{T_1},\ \forall j=1,\ldots,n_1\}\subset\Phi_1^{-1}({\mathcal K}^{T_1})\]
while $\{F\in{\mathcal K}^T:Y(u_{|V_j})=Y_j^{T_1},\ \forall j=1,\ldots,n_1\}$ is a nonempty open subset of ${\mathcal K}^T$. Therefore, $\Phi_1({\mathcal K}^T)\subset{\mathcal K}^{T_1}$. Similarly, we have
\[\{F\in{\mathcal K}^T:Y(u_{|V_j/V_{n_1}})=Y_{j-n_1}^{T_2},\ \forall j=n_1+1,\ldots,n_1+n_2\}\subset\Phi_2^{-1}({\mathcal K}^{T_2})\]
while $\{F\in{\mathcal K}^T:Y(u_{|V_j/V_{n_1}})=Y_{j-n_1}^{T_2},\ \forall j=n_1+1,\ldots,n_1+n_2\}$
is a nonempty open subset of ${\mathcal K}^T$. Therefore, $\Phi_2({\mathcal K}^T)\subset{\mathcal K}^{T_2}$.
It results $\Phi({\mathcal K}^T)\subset {\mathcal K}^{T_1}\times{\mathcal K}^{T_2}$.
The product ${\mathcal K}^{T_1}\times{\mathcal K}^{T_2}$ is an irreducible variety,
and due to formula (\ref{relation-dimension-Bu}), it has the same dimension as ${\mathcal K}^T$.
Hence, we actually have the equality $\Phi({\mathcal K}^T)={\mathcal K}^{T_1}\times{\mathcal K}^{T_2}$.
We have finally obtained that $\Phi$ restricts to an isomorphism between
${\mathcal K}^T$ and the product ${\mathcal K}^{T_1}\times{\mathcal K}^{T_2}$.
The proof of the proposition is then complete.
\hfill $\Box$

\subsection{Proof of Theorem \ref{theorem-3}}

The proof relies on Proposition \ref{proposition-somme-tableaux} and the following additional
preliminary result (see \cite[Proposition 6.1]{Fresse-Melnikov-2}).

\begin{lemma}
\label{lemma-inductive-fiber-bundle}
Let $T$ be a standard tableau of entries $1,\ldots,n$,
and let $T'$ be the subtableau of entries $1,\ldots,n-1$.
Let ${\mathcal K}^T,{\mathcal K}^{T'}$ be the components corresponding to $T,T'$
in the appropriate Springer fibers.
Assume that $n$ lies in the last column of $T$.
Let $k$ denote the length of the last column of $T$.
Then ${\mathcal K}^T$ is a locally trivial fiber bundle of base the projective space $\mathbb{P}^{k-1}$,
of fiber isomorphic to ${\mathcal K}^{T'}$.
\end{lemma}

\noindent{\em Proof of Lemma \ref{lemma-inductive-fiber-bundle}.}
Let $s$ be the nilpotent index of $u$ (which coincides with the number of columns of $T$).
We then have $k=\dim (V/\ker u^{s-1})$.
In the case where $n$ lies in the last column of $T$, it follows from the definition of ${\mathcal K}^T$
that $\ker u^{s-1}\subset V_{n-1}$ for all $F=(V_0,\ldots,V_n)\in{\mathcal K}^T$.
Then, we consider the map
\[\Phi:{\mathcal K}^T\rightarrow\mathbb{P}^{k-1},\ \ (V_0,\ldots,V_n)\mapsto V_{n-1}/\ker u^{s-1}\]
and we show that it is a locally trivial fiber bundle of fiber ${\mathcal K}^{T'}$
(see \cite[Proposition 6.1]{Fresse-Melnikov-2}, or \cite[Theorem 2.1]{Fresse-Melnikov-1}).
\hfill $\Box$

\bigskip

Now, let us prove Theorem \ref{theorem-3}.
We fix $V=\mathbb{C}^n$ and a nilpotent endomorphism $u\in\mathrm{End}(V)$.
Let $\lambda(u)=(\lambda_1,\ldots,\lambda_r)$ be the sizes of its Jordan blocks
and let $Y(u)$ be the corresponding Young diagram (i.e., of rows of sizes $\lambda_1,\ldots,\lambda_r$).
Let $\lambda^*(u)=(\lambda^*_1,\ldots,\lambda^*_s)$ be the conjugate partition of $n$.
Let $u^*\in\mathrm{End}(V)$ be nilpotent, such that $u,u^*$ have conjugate Jordan forms,
hence $\lambda(u^*)=\lambda^*(u)$.
Let $Y(u^*)$ be the Young diagram corresponding to $u^*$.
If $T$ is a standard tableau of shape $Y(u)$, then its transposed tableau $T^*$ has shape $Y(u^*)$.
Hence, to the irreducible component ${\mathcal K}^T\subset{\mathcal B}_u$,
we associate the component ${\mathcal K}^{T^*}\subset{\mathcal B}_{u^*}$.
We suppose that ${\mathcal K}^T$ contains a dense Jordan orbit.
Our purpose is to show that ${\mathcal K}^{T^*}$ is an iterated fiber bundle
of base $(\mathbb{P}^1,\ldots,\mathbb{P}^{\lambda_1-1},\ldots,\mathbb{P}^1,\ldots,\mathbb{P}^{\lambda_r-1})$.
We reason by induction on $n$, with immediate initialization for $n=1$.
Suppose the property holds until the rank $n-1\geq 1$ and let us show for $n$.

By Proposition \ref{proposition-2}, the tableau $T$ is of the form $T=T_\pi$ for some $\pi\in\Pi^1_u$.
Set $\pi=(I_1,\ldots,I_r)$, where $|I_j|=\lambda_j$ and $I_1\sqcup\ldots\sqcup I_r=\{1,\ldots,n\}$.
We distinguish two cases, depending on whether $1,n$ belong to the same subset $I_j$ or not.

\medskip
\noindent (1) Suppose there is $j\in\{1,\ldots,r\}$ such that $1,n\in I_j$.

\smallskip

Then, we have $I_k<I_j$ for all $k\not=j$ (see section \ref{generalized-bc-components}).
By definition of the set $\Pi_u^1$, this implies that $|I_j|=\min\{|I_k|:k=1,\ldots,r\}$.
In other words, $\lambda_j$ is the length of the minimal row of the tableau $T_\pi$.
By definition of $T_\pi$ (see section \ref{Jordan-orbits}), $n$ is the last entry of the $\lambda_j$-th column
of $T_\pi$, hence it lies in the last row of $T_\pi$.
By consequent, $n$ lies in the last column of the transposed tableau $T^*=(T_\pi)^*$.

Note that the last column of $T^*$ has length $\lambda_r$.
Let $(T^*)'$ be the subtableau of $T^*$ of entries $1,\ldots,n-1$. By Lemma \ref{lemma-inductive-fiber-bundle},
there is a fiber bundle ${\mathcal K}^{T^*}\rightarrow \mathbb{P}^{\lambda_r-1}$ of fiber
isomorphic to ${\mathcal K}^{(T^*)'}$.

Let $\lambda'_r=\lambda_r-1$ and for $k\not=r$, let $\lambda'_k=\lambda_k$.
Similarly as in section \ref{section-extremal-point},
we consider $u'\in\mathrm{End}(\mathbb{C}^{n-1})$ of Jordan form $\lambda(u')=(\lambda'_1,\ldots,\lambda'_r)$.
Set $I'_j=I_j\setminus\{n\}$, and $I'_k=I_k$ for $k\not=j$.
Let $\pi'\in(I'_1,\ldots,I'_r)$. As in section \ref{section-extremal-point}, we see that $\pi'\in\Pi_{u'}^1$.
Moreover, the tableau $T_{\pi'}$ is the subtableau of $T_\pi$ of entries $1,\ldots,n-1$.
It follows that $(T^*)'$ coincides the transposed tableau $(T_{\pi'})^*$.
Then, by induction hypothesis, the component ${\mathcal K}^{(T^*)'}$ is an iterated fiber bundle
of base $(\mathbb{P}^1,\ldots,\mathbb{P}^{\lambda_1-1},\ldots,\mathbb{P}^1,\ldots,\mathbb{P}^{\lambda_r-2})$.
Therefore, ${\mathcal K}^{T^*}$ is an iterated fiber bundle of the desired base.

\medskip

\smallskip
\noindent (2) Suppose there are $j,k\in\{1,\ldots,r\}$, $j\not=k$ such that $1\in I_j$, $n\in I_k$.

\smallskip

We set in this case $n_1=\max I_j$. As $n\notin I_j$, we have $n_1<n$.
Again according to the definition of the set
$\Pi_u^1$, the sequence $\pi$ has no crossing, which implies that for all $l$, we have
either $I_l\subset [1,n_1]$ or $I_l\subset[n_1+1,n]$.
We consider the two subsequences $\pi_1=(I_l:l=1,\ldots,r,\ I_l\subset [1,n_1])$
and $\pi_2=(\phi(I_l):l=1,\ldots,r,\ I_l\subset [n_1+1,n])$,
where $\phi:\{n_1+1,\ldots,n\}\rightarrow\{1,\ldots,n-n_1\}$ is defined by $\phi(i)=i-n_1$.
Let $u_1\in\mathrm{End}(\mathbb{C}^{n_1})$ and $u_2\in\mathrm{End}(\mathbb{C}^{n-n_2})$ be
nilpotent of Jordan forms $\lambda(u_1)=(\lambda_l:l=1,\ldots,r,\ I_l\subset [1,n_1])$
and $\lambda(u_2)=(\lambda_l:l=1,\ldots,r,\ I_l\subset [n_1+1,n])$ respectively.
Then, we clearly have $\pi_1\in\Pi_{u_1}^1$ and $\pi_2\in\Pi_{u_2}^1$.

Let $T_1=T_{\pi_1}$ and $T_2=T_{\pi_2}$ be the standard tableaux associated to $\pi_1$ and $\pi_2$,
and let $(T_1)^*$ and $(T_2)^*$ be their respective transposes.
From the definition of the tableaux $T_\pi,T_{\pi_1},T_{\pi_2}$, it is easy to see
that the transpose $T^*=(T_\pi)^*$ is in fact the sum $T^*=(T_1)^*+(T_2)^*$.
By Proposition \ref{proposition-somme-tableaux},
it follows that the component ${\mathcal K}^{T^*}$
is isomorphic to the product ${\mathcal K}^{(T_1)^*}\times {\mathcal K}^{(T_2)^*}$.
Moreover, by induction hypothesis, ${\mathcal K}^{(T_1)^*}$ is an iterated fiber bundle
of base $(\mathbb{P}^1,\ldots,\mathbb{P}^{\lambda_l-1})_{l}$ where $l$ runs over
$\{l=1,\ldots,r:I_l\subset[1,n_1]\}$. Similarly
${\mathcal K}^{(T_2)^*}$ is an iterated fiber bundle
of base $(\mathbb{P}^1,\ldots,\mathbb{P}^{\lambda_l-1})_{l}$ for $l$ running over
$\{l=1,\ldots,r:I_l\subset[n_1+1,n]\}$. We infer that ${\mathcal K}^{T^*}$ is an iterated fiber
bundle of the desired base.

\smallskip

The proof of Theorem \ref{theorem-3} is then complete.

\section*{Index of the notation}

\begin{tabbing}
\ \ \= \S\ref{outline}\qquad \= $V$, $n$, $u$, ${\mathcal B}$, ${\mathcal B}_u$, \\[1mm]
\> \ref{section-Yu} \> $\lambda(u)$, $Y(u)$, $\lambda^*(u)$, \\[1mm]
\> \ref{Spaltenstein-construction} \> $Y_i^T$, $Y(u_{|V_i})$, ${\mathcal B}_u^T$, ${\mathcal K}^T$, \\[1mm]
\> \ref{def-BC-R} \> $\Lambda_u$, $\Lambda^*_u$, ${\mathcal U}_\pi^{\mathrm{BC}}$, ${\mathcal K}_\pi^{\mathrm{BC}}$, ${\mathcal K}_\pi^{\mathrm{R}}$, \\[1mm]
\> \ref{duality} \> $T_\pi$ ($\pi\in\Lambda_u$), $T_\pi^*$ ($\pi\in\Lambda_u^*$), $u^*$, $T^*$, \\[1mm]
\> \ref{first-result} \> $\pi\geq\rho$, \\[1mm]
\> \ref{second-result} \> $Z_u$, \\[1mm]
\> \ref{Jordan-orbits} \> $\Pi_u$, $\pi:\{1,\ldots,n\}\rightarrow\{\emptyset,1,\ldots,n\}$, ${\mathcal Z}_\pi$, $T_\pi$ ($\pi\in\Pi_u$), $\Pi_u^0$, \\[1mm]
\> \ref{generalized-bc-components} \> $I_j<I_k$, $\Pi_u^1$, \\[1mm]
\> \ref{flags-Fd} \> $\tau_0$, $\sigma(\tau_0)$, $\mathbf{S}_u$, $F_\sigma$, $D$,
$\tau_{\underline{d}}$, $\sigma_{\underline{d}}$, $F_{\underline{d}}$, \\[1mm]
\> \ref{proof-proposition-p-1-p-1} \> $\Omega_{\underline{d}}$, $\phi_{i,j}$, \\[1mm]
\> \ref{sum-tableaux} \> $Y_1+Y_2$, $T_1+T_2$.
\end{tabbing}

\end{document}